\documentclass[12pt]{nmeauthNoIJNME} 

\usepackage[colorlinks, linktocpage=true, pdfpagemode=No,
linkcolor=blue, citecolor=blue, pagecolor=blue]{hyperref}

\usepackage{listings}
\usepackage{multirow}
\usepackage{csml}
\usepackage{float}
\usepackage{subfig}
\usepackage{booktabs}
\usepackage[table]{xcolor}
\usepackage{mathtools}
\newcommand{\RVE}{\operatorname{\mathit{R\kern-.17em V\kern-.17em E}}}
\begin{document}

\runningheads{S. Herath, X. Xiao and F. Cirak} {Data-driven homogenisation of knitted membranes}


\title{Computational modelling and data-driven homogenisation of knitted membranes} 

\author{Sumudu Herath\affil{1,3}, Xiao Xiao\affil{2,3}\corrauth and Fehmi Cirak\affil{3}}

\address
{\affilnum{1}Department of Civil Engineering, University of Moratuwa, Moratuwa, Sri Lanka \break
 \affilnum{2}Inria, 2004 route des Lucioles, 06902 Sophia Antipolis, France \break
 \affilnum{3}Department of Engineering, University of Cambridge, Trumpington Street, Cambridge CB2 1PZ, U.K.
}

\corraddr{Inria, 2004 route des Lucioles, 06902 Sophia Antipolis, France. E-mail: xiao.xiao@inria.fr}

\begin{abstract}		
Knitting is an effective technique for producing complex three-dimensional surfaces owing to the inherent flexibility of interlooped yarns and recent advances in manufacturing providing better control of local stitch patterns. Fully yarn-level modelling of large-scale knitted membranes is not feasible. Therefore, we use a two-scale homogenisation approach and model the membrane as a Kirchhoff-Love shell on the macroscale and as Euler-Bernoulli rods on the microscale. The governing equations for both the shell and the rod are discretised with cubic B-spline basis functions. For homogenisation we consider only the in-plane response of the membrane. The solution of the nonlinear microscale problem requires a significant amount of time due to the large deformations and the enforcement of contact constraints, rendering conventional online computational homogenisation approaches infeasible. To sidestep this problem, we use a pre-trained statistical Gaussian Process Regression (GPR) model to map the macroscale deformations to macroscale stresses. During the offline learning phase, the GPR model is trained by solving the microscale problem for a sufficiently rich set of deformation states obtained by either uniform or Sobol sampling. The trained GPR model encodes the nonlinearities and anisotropies present in the microscale and serves as a material model for the membrane response of the macroscale shell. The bending response can be chosen in dependence of the mesh size to penalise the fine out-of-plane wrinkling of the membrane. After verifying and validating the different components of the proposed approach, we introduce several examples involving membranes subjected to tension and shear to demonstrate its versatility and good performance.
\end{abstract}


\maketitle

%

%
\section{Introduction}
\label{sec:introduction}

Knitting is one of the most efficient and widely used techniques for producing fabric membranes. The recent advances in computational knitting make it possible to produce large complex three-dimensional surfaces in one piece without seams~\cite{mccann2016compiler,popescu2018automated}. On a modern programmable flat-bed knitting machine it is possible to produce even non-developable surfaces by a local variation of the stitch pattern consisting of an interlooped yarn. Most promisingly, the yarn can be replaced or integrated with electroactive or conductive yarns to produce novel interactive textiles with sensing and/or actuation capabilities~\cite{luo2021knitui,maziz2017knitting,persson2018actuating}. Knitted membranes are usually very flexible because the primary deformation mechanism for the yarn is bending rather than axial stretching. Their unique stretchability and drapability properties make knitted membranes appealing as a reinforcement in composite components~\cite{leong2000potential,Hasani2017}. If needed, the stiffness can be increased by inserting straight high-strength fibres during the knitting process. Such reinforced membranes have been recently used as a formwork in architectural engineering~\cite{popescu2018building}. Considering the recent advances in knitting, there is a need for efficient computational approaches for the analysis of large-scale knitted membranes. 

For large-scale analysis of knitted membranes, the computational homogenisation approaches which take into account the deformation of the interlooped yarn on the microscale are crucial. There is an extensive amount of literature on finite element-based computational homogenisation of heterogeneous solids, see e.g. the reviews~\cite{peric2011micro,geers2017,saeb2016aspects}. Two-scale homogenisation often referred to as FE$^\text{2}$, combined with a yarn-level and a membrane-level finite element model has been also applied to woven and knitted fabrics~\cite{nadler2006,fillep2013computational,liu2019multiscale,do2020nonlinear}. The boundary conditions of the microscale representative volume element (RVE) are given by the membrane deformation and in turn the averaged yarn stress in the RVE yields the membrane stress. However, such schemes are inefficient for knitted membranes with large deformations because of the need to solve a nonlinear problem at each quadrature point of the membrane. It is increasingly apparent that for nonlinear problems two-scale homogenisation must be considered in combination with a data-driven machine learning model~\cite{bessa2017}. The model can be trained in an off-line learning phase by solving the microscale problem for a sufficiently rich set of deformation states. Subsequently, the trained model provides a closed-form constitutive equation that is used in the macroscale model. As a machine learning model, for instance, neural networks, GPR or polynomial regression have been used~\cite{bessa2017,le2015computational,nguyen2020surrogate,wu2020recurrent,do2020nonlinear}. Alternatively, it is possible to formulate the  macroscale finite element problem directly on the training data set bypassing the need to train a machine learning model~\cite{karapiperis2021data,platzer2021finite}.
GPR is a very widely used and well-studied Bayesian statistical technique and provides as such a consistent approach for dealing with epistemic and aleatoric uncertainties and issues such as overfitting~\cite{rasmussen2006}. In contrast to most other regression techniques, GPR provides an estimate of the uncertainty in the prediction and a principled approach to learn any model hyperparameters through evidence maximisation. Consequently, GPR is well suited for highly nonlinear, multi-modal and multi-dimensional regression problems from engineering; see~\cite{forrester2008engineering} for a comparison of different regression techniques. Therefore, we chose to approximate the response of the microscale yarn model with GPR. Their limitation to relatively small dimensions and number of training points is not relevant for this paper.

Owing to the relative slenderness of the yarn on the microscale and the membrane on the macroscale, they are best modelled as a rod and a shell, respectively. Evidently, this leads in comparison to a 3D solid model to an immense reduction in the number of degrees of freedom. The geometrically exact rod theories pioneered by Simo et al.~\cite{simo1985,simo1986} provide a consistent and efficient framework for modelling rods undergoing finite deformations. Later contributions to geometrically exact rod theories include~\cite{cardona1988,KONDOH1986253,jelenic1999geometrically}. While in most of these classical rod models the  transverse shear deformations are taken into account, in more recent Euler-Bernoulli type models they are neglected~\cite{meier2014,bauer2016nonlinear}. The omission of the transverse shear effectively sidesteps the shear locking problem which is a major impediment in the analysis of slender rods. These new models are usually discretised using smooth B-spline basis functions due to the presence of the higher-order displacement derivatives in the energy functional. The Euler-Bernoulli type rod model introduced in this paper takes into account the stretching, bending and torsion of the yarn as well as the non-frictional rod-to-rod contact between yarns. Similar to the yarn, the fabric membrane can be modelled with geometrically exact shell theories going back again to Simo et al.~\cite{Simo:1989aa,Simo:1989ab}. As for rods, in more recent Kirchhoff-Love type shell models the transverse shear deformations are neglected and the weak form is discretised using smooth B-splines or subdivision basis functions~\cite{cirak2000,Cirak:2001aa,Kiendl:2009aa}.  In this paper we make use of the Kirchhoff-Love subdivision shell implementation introduced in~\cite{Cirak:2011aa,long2012}. The true bending stiffness of knitted membranes is usually very small so that under compression excessively small wrinkles may form. The need to resolve a large number of small wrinkles leads in engineering computations often to an impractically large number of elements. A sufficiently large (artificial) bending stiffness can act as an regulariser and render the problem computable by increasing the size of the wrinkles~\cite{Cirak:2001aa,cirak2014computational}. Alternatively, using a tension field approach the wrinkling of the membrane may be suppressed by omitting compressive stresses so that element sizes can be chosen reasonably large; see~\cite{wong2006wrinkled2,wong2006wrinkled} for a discussion on different approaches in the case of non-knitted membranes.

The outline of the paper is as follows. In Section~\ref{sec:3Drods} we introduce our finite deformation rod model, its discretisation with B-splines as well as the treatment of rod-to-rod contact. Subsequently, in Section~\ref{sec:compHom} we discuss the proposed microscale yarn-level RVE model for computational homogenisation. This model is verified and validated with experimental and numerical results from literature. In Section~\ref{sec:knitGPRhomogen}, we introduce the data-driven GPR model and describe its training with the yarn-level RVE model. Finally, in Section~\ref{sec:examples}, we first discuss the training of the GPR model and then analyse membranes subjected to tension to assess the accuracy of the obtained data-driven constitutive model.

%

\section{Finite deformation analysis of yarns} \label{sec:3Drods}
%
In this section we summarise the governing equations for the finite deformation Euler-Bernoulli rod model for the yarn. The presented equations are without loss of generality restricted to rods with circular cross-sections. We take into account rod-to-rod contact by enforcing the non-penetration constant with the Lagrange multiplier method. 

\subsection{Kinematics} \label{subsec:rodkinematics}
%
The geometry of the rod is described by a set of circular cross-sections connected by their line of centroids. In accordance with the Euler-Bernoulli assumption, equivalent to the Kirchhoff-Love assumption for shells and plates, transverse shear is neglected so that cross-sections remain always normal to the line of centroids.

\begin{figure}
	\centering
  \includegraphics[width=0.85\columnwidth]{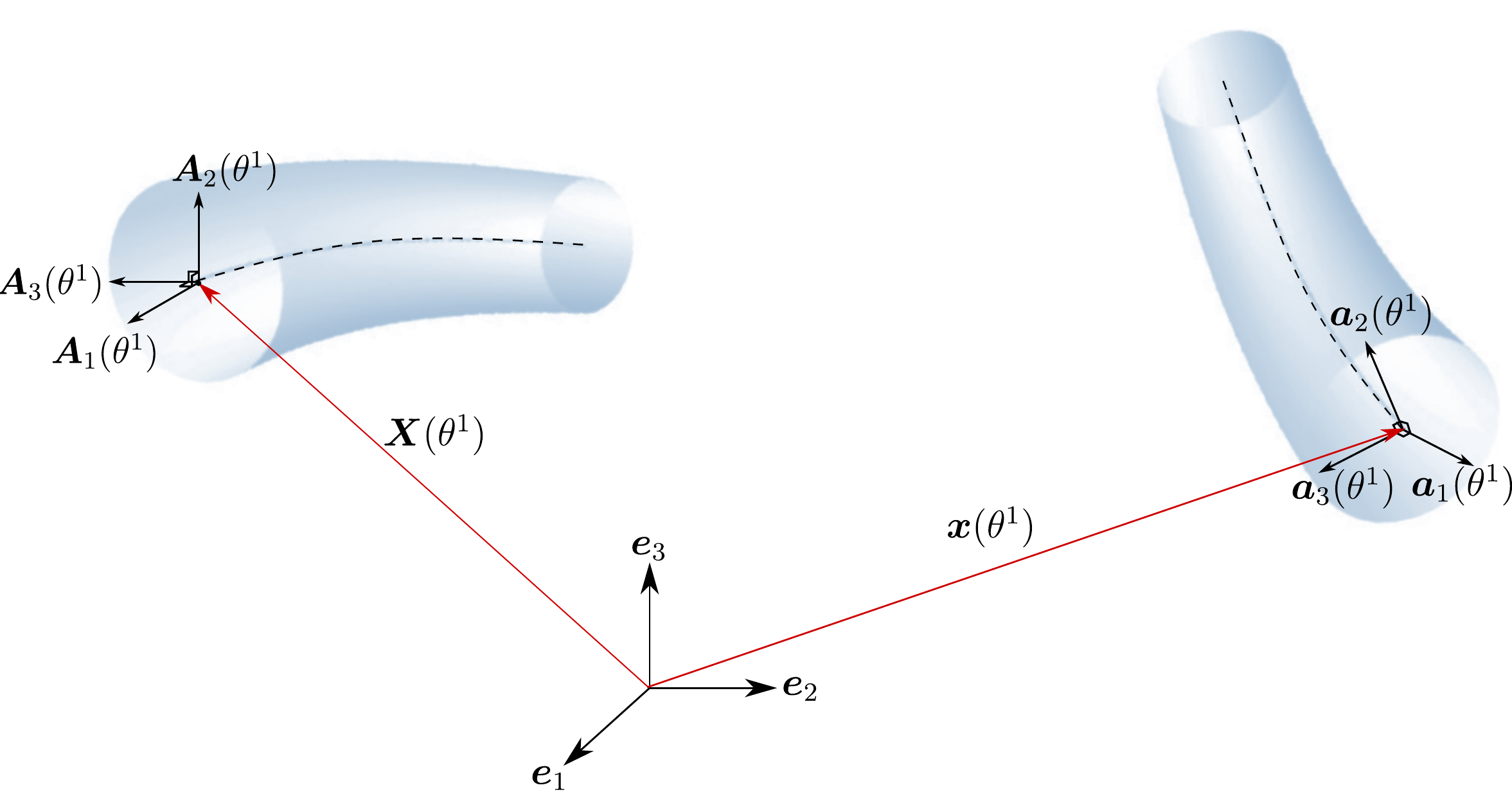}
  \caption{Geometric description of a spatial rod in its reference (left) and deformed (right) configurations. The two configurations are defined using the line of centroids~$\vec X(\theta^1)$ and~$\vec x(\theta^1)$ and the respective covariant basis vectors.}
  \label{fig:spatialrod}
\end{figure}

The position vectors of material points in the reference and deformed configurations \mbox{$\vec{R}(\theta^1, \, \theta^2, \, \theta^3)$} and \mbox{$\vec{r} (\theta^1, \, \theta^2, \, \theta^3) $} are parametrised in terms of the convective coordinates~\mbox{$\{\theta^1, \, \theta^2, \, \theta^3 \}$} as
\begin{subequations}  \label{eq:rodKinematics}
\begin{align}
	\vec{R}(\theta^1, \, \theta^2, \, \theta^3) &= \vec{X}(\theta^1) + \theta^2 \vec{A}_2(\theta^1)  + \theta^3 {\vec{A}}_3(\theta^1) \label{eq:rodMaterialPointsRef} \, , \\ 
	\vec{r}(\theta^1,  \, \theta^2, \,  \theta^3) &= \vec{x}(\theta^1) + \theta ^2 \vec{a}_2(\theta^1)  + \theta^3 {\vec{a}}_3(\theta^1) \label{eq:rodMaterialPointsCur} \, , 
\end{align}
\end{subequations} 
where~$\vec{X}(\theta^1)$ and~$\vec{x}(\theta^1) $  denote the lines of centroids parameterised by~$\theta^1$,  see Figure \ref{fig:spatialrod}. In turn, the cross-section  with the radius~$r$  is parameterised by~\mbox{$| \theta^2 | \le r$} and \mbox{$ |  \theta^3 | \le r $}. The tangent vectors to the line of centroids are given by
\begin{equation}
	\vec A_1 = \frac{\text{d} \vec X}{\text{d} \theta^1} = \vec{X}_{,1}  \, ,  \quad  \vec a_1 = \frac{\text{d} \vec x}{\text{d} \theta^1} = \vec{x}_{,1}  \, . \label{eq:a2a3MApping}
\end{equation}
The two orthonormal directors~$\vec A_2$ and~$\vec A_3$ are chosen so that they satisfy in the reference configuration
\begin{equation}
	| \vec A _\iota | =1 \, ,  \quad  \vec A_2  \cdot \vec A_3 = 0 \, , \quad \vec A_1 \cdot \vec A_\iota  = 0 \, .
\end{equation}
Here and in the following Greek indices take the values~$ \{ 2, \, 3 \}$ and summation over repeated indices is assumed. The two orthonormal directors~$\vec a_\iota$ in the deformed configuration are obtained by rotating the reference configuration directors with a rotation matrix~$\vec \Lambda \in \text{SO}_3$, 
\begin{equation}
	\vec{a}_{\iota} = \vec \Lambda \vec{A}_{\iota}  \, .
\end{equation}
This, in combination with the Euler-Bernoulli assumption, ensures that the two directors~$\vec a_2$ and~$\vec a_3$ in the deformed configuration satisfy
\begin{equation}
	| \vec a _\iota | =1 \, ,  \quad  \vec a_2  \cdot \vec a_3 = 0 \, , \quad \vec a_1 \cdot \vec a_\iota  = 0 \ .
\end{equation}
The rotation $ \vec \Lambda$ is composed of two rotations,
\begin{equation} \label{eq:rotationMatrix}
	 \vec \Lambda (\vec x, \, \vartheta ) = \vec \Lambda_2 (\vec x) \vec \Lambda_1 (\vartheta) \, .
\end{equation}
That is, the reference directors $\vec A_\iota$ are mapped to the deformed directors $\vec a_\iota$ in two steps using an intermediary configuration with $\vec a'_\iota = \vec \Lambda_1( \vartheta ) \vec A_\iota$. The matrix $\vec \Lambda_1 (\vartheta)$ describes a rotation by an angle~$\vartheta$ about the unit tangent $\hat {\vec A}_1 = \vec A_1 / |\vec A_1| $ and is according to the Rodrigues formula given by 
\begin{equation} \label{eq:R1} 
	 \vec \Lambda_1(\vartheta ) = \vec{I} + \sin\vartheta \vec{L} + (1-\cos\vartheta) \vec{L}\vec{L} 
\end{equation}
with the identity matrix $\vec I$ and the skew-symmetric matrix 
\begin{equation}	
	\vec{L}= 
	\begin{pmatrix*}[c]
	0 & -\hat {\vec{A}}_{1} \cdot \vec e_3 & \phantom{+} \hat {\vec{A}}_{1} \cdot \vec e_2 \\
	\phantom{+} \hat {\vec{A}}_{1}  \cdot \vec e_3 & 0 & -\hat {\vec{A}}_{1} \cdot \vec e_1& \\
	-\hat {\vec{A}}_{1} \cdot \vec e_2 &\phantom{+} \hat {\vec{A}}_{1} \cdot \vec e_1 & 0 \\
\end{pmatrix*}
\, .
\end{equation}
Subsequently, we use the smallest rotation formula~\cite{meier2014, crisfield1997non} for the second matrix~$\vec \Lambda_2 (\vec a_1)$, which maps~$\hat {\vec A}_1$ to the unit tangent vector~$\hat {\vec a}_1 = \vec a_1 / |\vec a_1| $,
\begin{align} \label{eq:R2}
	\vec \Lambda_2(\vec a_{1} ) &= \vec{I} - \dfrac{(\hat{\vec A}_{1}+\hat{\vec a}_{1}) \otimes \hat{\vec a}_{1}}{1+\hat{\vec a}_{1} \cdot \hat{\vec A}_{1}} \, . 
\end{align}

To derive the strains corresponding to the assumed kinematics~\eqref{eq:rodKinematics}, we consider the Green-Lagrange strain tensor of 3D elasticity 
\begin{align}\label{eq:GreenLagStrain}
	\vec{E} = E_{ij }  \  {\vec{G}}^i\otimes{\vec{G}}^j      = \frac{1}{2}(g_{ij}-{G}_{ij}) \ {\vec{G}}^i\otimes{\vec{G}}^j  \, .
\end{align}
Here and in the following Latin indices take the values~$\{ 1, \, 2 , \, 3 \}$. The covariant basis vectors~$\vec G_i$ and~$\vec g_i$  and the contravariant basis vectors~$\vec G^i$ and $\vec g^i$ are defined as 
\begin{align}\label{eq:covariantBasevectors}
	{\vec{G}}_i &=  \frac{\partial {\vec R}}{\partial \theta^i} \, , \quad  \vec{g}_i =  \frac{\partial {\vec r}}{\partial \theta^i} \,,  \quad  \vec G_i \cdot {\vec{G}}^j =  \delta_i^j \, , \quad \vec g_i \cdot  {\vec{g}}^j =  \delta_i^j ,
\end{align}
where~$\delta_i^j$ is the Kronecker delta. The corresponding two metric tensors~${G}_{ij}$ and~$g_{ij}$ are given by
\begin{equation}
	G_{ij} = \vec G_i \cdot  \vec G_j \, ,  \quad  g_{ij} = \vec g_i \cdot \vec g_j \, .
\end{equation}
After introducing the assumed kinematics~\eqref{eq:rodKinematics} in~\eqref{eq:GreenLagStrain} and some algebraic simplifications we obtain for the components of the strain tensor 
\begin{subequations}\label{eq:allStrainComponents}
\begin{align}
	[E_{ij}] &= 
	\begin{bmatrix}
		\alpha + \theta^2 \beta_{2} + \theta^3 \beta_{3} & \theta^3 \gamma  & -\theta^2 \gamma \\
		 & 0 & 0 \\
		\text{symmetric} &  & 0 \\
	\end{bmatrix}
\intertext{with} 
	\alpha(\vec x)& = \frac{1}{2} (\vec{a}_1 \cdot \vec{a}_1 - {\vec{A}}_1 \cdot {\vec{A}}_1) \, , \\
	\beta_{2}(\vec x, \vartheta)& = {\vec{A}}_2 \cdot {\vec{A}}_{1,1} - \vec{a}_2 \cdot \vec{a}_{1,1} \, , \\
	\beta_{3}(\vec x, \vartheta)& = {\vec{A}}_3 \cdot {\vec{A}}_{1,1} - \vec{a}_3 \cdot \vec{a}_{1,1} \, , \\
	\gamma(\vec x, \vartheta) & = \frac{1}{2}(\vec a_2 \cdot \vec a_{3,1} - \vec A_2 \cdot \vec A_{3,1} ) \, .
\end{align}
\end{subequations}
We identify $\alpha, \ \beta_{2}, \ \beta_{3} \ \text{and} \ \gamma$ as the membrane, bending about $\vec a_2$, bending about $\vec a_3$ and torsional shear strains, respectively. Moreover, due to the Euler-Bernoulli assumption the strain components $E_{23}, \ E_{32}$ and $E_{33}$ are zero and the in-plane shear strains $E_{12}$ and $E_{13}$ are induced only by torsional shear strain~$\gamma$.

%
\subsection{Equilibrium equations in weak form}\label{subsec:rod_energyDiscretisation}
%
The potential energy of a rod with the line of centroids~\mbox{$\Gamma \subset \mathbb R$} and the cross-section~\mbox{$\Omega \subset \mathbb R^2$} occupying the volume~\mbox{$V = \Gamma \times \Omega \subset \mathbb R^3$} in its reference configuration takes the form 
\begin{align} \label{eq:potentialEnergy}
	\Pi(\vec r) = \Pi_{\text{int}}(\vec r) + \Pi_{\text{ext}}(\vec r)  = \int_{V} \psi ( \vec E ) \D V + \Pi_{\text{ext}}(\vec r)  \ , 
\end{align}
where $\psi(\vec E)$ is the strain energy density and~$ \Pi_{\text{ext}}(\vec r)$ is the potential of the externally applied forces. At equilibrium the potential energy of the rod is stationary, i.e.,  
\begin{subequations} \label{eq:variation}
\begin{align} 
	\delta \Pi(\vec r) &=  \delta \Pi_{\text{int}}(\vec r) + \delta \Pi_{\text{ext}}(\vec r)  = 0    \label{eq:totVariation}
\intertext{with the external virtual work~$\delta \Pi_{\text{ext}}(\vec r)$ and the internal virtual work}
	\delta \Pi_{\text{int}}(\vec r) &= \int_V \frac{ \partial \psi( \vec E) }{ \partial \vec E } :  \delta \vec E  \D V   =  \int_V \vec S  : \delta \vec E  \D V   \, ,  \label{eq:intVariation}
\end{align}
\end{subequations}
where $\vec{S}$ is the second Piola-Kirchhoff stress tensor. The strain tensor of the rod~\eqref{eq:allStrainComponents} depends on the displacement of the line of centroids 
\begin{equation} \label{eq:centroidDisp}
	\vec u (\theta^1) = \vec x  (\theta^1) - \vec X  (\theta^1) \, ,
\end{equation}
and the rotation angle~$\vartheta$. Hence, we can write for the internal virtual work~\eqref{eq:intVariation} more succinctly 
\begin{subequations} \label{eq:DELTAinernalEnergy1}
\begin{align}
	\delta \Pi_{\text{int}}(\vec x) &=  \int_{\Gamma} \int_{\Omega}  \vec S : \dfrac{\partial \vec E}{\partial \vec u} \cdot \delta \vec u | \vec A_1 |  \D \Omega \D \Gamma + \int_{\Gamma}\int_{\Omega} \vec S : \dfrac{\partial \vec E}{\partial \vartheta}  \delta \vartheta | \vec A_1 |  \D \Omega \D \Gamma \\ &= 
	\int_\Gamma \vec f_{\vec u} \cdot \delta \vec u | \vec A_1 |  \D \Gamma + \int_\Gamma f_{\vartheta}  \delta \vartheta   | \vec A_1 |  \D \Gamma \, ,
\end{align}
\end{subequations}
where $\vec f_{\vec u}$ and $f_{\vartheta}$ are the internal forces conjugate to the virtual displacements~$\delta \vec u$ and rotations~$\delta \vartheta$.

As a material model we use the isotropic St Venant-Kirchhoff model with the strain energy density 
\begin{align}  \label{eq:stVenantEnergyDens}
	\psi(\vec E) = \dfrac{1}{2} \vec E : \vec C : \vec E = \dfrac{1}{2} \vec S : \vec E \, , 
\end{align}
and the fourth-order constitutive tensor
\begin{equation} \label{eq:stVKb}
\vec C = C^{ijkl} \ \vec G_i \otimes \vec G_j \otimes  \vec G_k  \otimes  \vec G_l  = \left[\lambda G^{ij} G^{kl} + \mu(G^{ik}G^{jl} + G^{il}G^{jk})\right] \vec G_i \otimes \vec G_j \otimes  \vec G_k  \otimes  \vec G_l \, ,
\end{equation}
where $\lambda$ and $\mu$ are the two Lam\'e parameters~\cite{ciarlet2006introduction}. The contravariant metric tensor~$G^{ij}$ is determined from the relation~$G_{ki}G^{ij} = \delta_k^j$. 

To derive analytical expressions for the internal forces we first introduce the rod strain tensor~\eqref{eq:allStrainComponents} and the constitutive equation~\eqref{eq:stVKb} in the internal virtual work~\eqref{eq:DELTAinernalEnergy1}. Subsequently, we integrate over the rod cross-section analytically to obtain the internal forces
\begin{subequations} \label{eq:internalForces}
\begin{align} 
	\vec f_{\vec u} & =  n\frac{\partial \alpha }{\partial \vec{u}}
	+ m_2 \frac{\partial \beta_{2} }{\partial \vec{u}} + m_3\frac{\partial \beta_{3} }{\partial \vec{u}} + q\frac{\partial \gamma }{\partial \vec{u}}  \, , \label{eq:internalforceU} \\
	f_{\vartheta} &=  m_{2}\frac{\partial \beta_{2} }{\partial \vartheta} + m_{3}\frac{\partial \beta_{3} }{\partial \vartheta} + q\frac{\partial \gamma }{\partial \vartheta} \, . \label{eq:internalforceTheta}
\end{align}
\end{subequations}
Here, the axial force $n$, the bending moments $m_\iota$ and the torque $q$ are defined as
\begin{subequations}
\begin{align}
 n&= E {G}^{11}{G}^{11} A \alpha  \, , \label{eq:membraneResultant} \\
 m_{\iota}&= E {G}^{11}{G}^{11} I_{\iota} \beta_{\iota} \, , \label{eq:bendingResultant}  \\
 q &= 2E {G}^{11} J  \gamma \, , \label{eq:torsionResultant} 
\end{align}
\end{subequations}
where $A= \pi r^2$, $I_\iota = \pi r^4/4$ and $J=\pi r^4/2$ are the cross-section area, second moment of area and torsional constant of the circular rod, and $E$ is its Young's modulus. The tedious but straightforward derivation of internal forces and their derivatives are summarised in Appendix~\ref{app:intForces}.
%
\subsection{Finite element discretisation}\label{subsec:FEMdiscretisation}
We follow the isogeometric analysis paradigm and  use univariate cubic B-splines to discretise the lines of centroids $ \vec X (\theta^1)$ and $\vec x (\theta^1)$ in the reference and deformed configurations. We choose smooth B-splines because the bending strains $\beta_\iota$ in~\eqref{eq:allStrainComponents} require at least $C^1$ continuous smooth basis functions. The kinematic relationship~\eqref{eq:centroidDisp} is restated after discretisation as
\begin{align}
\sum_{I=1}^{n_B} B^I(\theta^1) \vec x_I = \sum_{I=1}^{n_B} B^I(\theta^1) \vec X_I  + \sum_{I=1}^{n_B} B^I(\theta^1) \vec u_I \, ,
\end{align}
where the B-spline basis $B^I(\theta^1)$ and their coefficients correspond to the $n_B$ control vertices on the discretised rod centreline. 

The discretised weak form of the equilibrium equations~\eqref{eq:variation} yields after linearisation an algebraic system of equations, which we solve with the Newton-Raphson scheme. For linearisation the gradients of the internal forces~\eqref{eq:internalForces} are required, see Appendix~\ref{App:HessianComp}.
%
\subsection{Yarn-to-yarn contact}\label{subsec:boundaryenforcement}
%
\begin{figure} 
\includegraphics[scale=0.7]{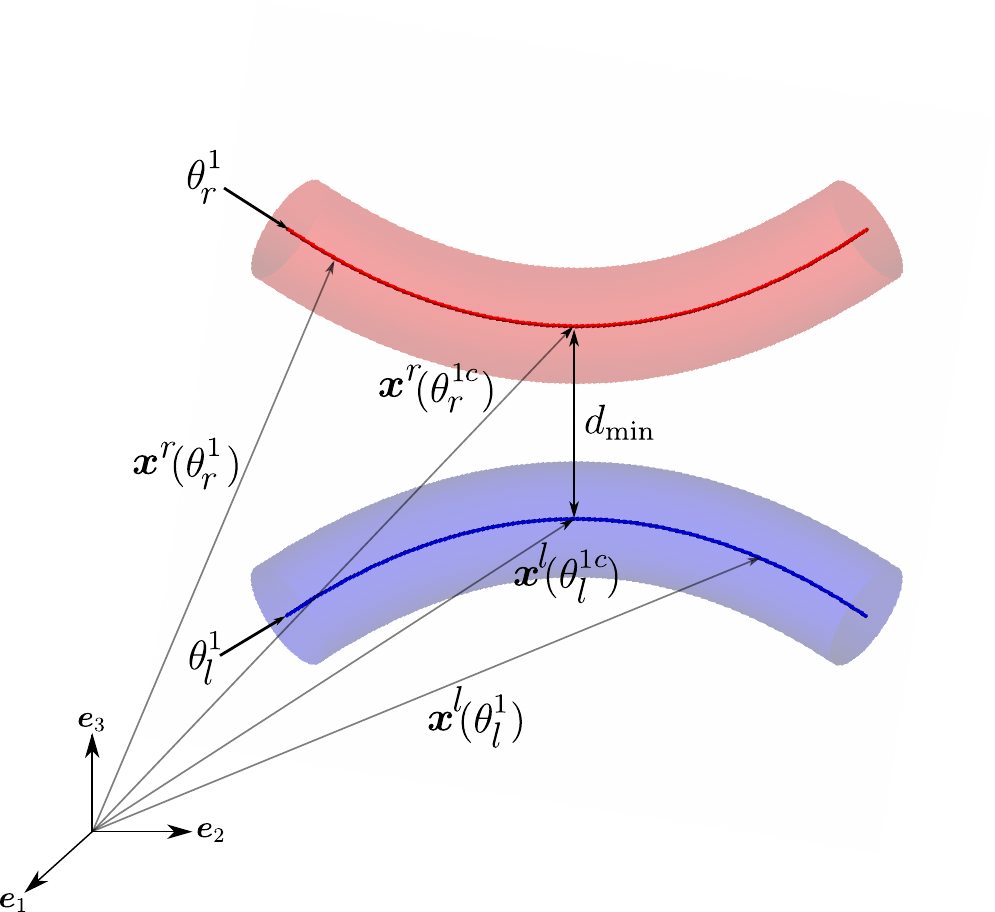}
\centering
\caption{Contact between two rods with the lines of centroids~$\vec x^l (\theta_l^1)$ and~$\vec x^r (\theta_r^1)$. Between the points~$\vec x^l (\theta_l^{1c})$  and~$\vec x^l (\theta_l^{1c})$ the distance is minimum and has the value~$d_{\min}$.  }
\label{fig:rodContact_P2P}
\end{figure}
We use the Lagrange multiplier method to consider pointwise non-frictional contact between two circular rods. By closely following Wriggers et al.~\cite{wriggers1997} and Weeger et al.~\cite{weeger2017}, we add to the total potential energy $\Pi(\vec r)$ in~\eqref{eq:potentialEnergy} the contact potential energy
\begin{equation}\label{eq:contact}
\Pi_{c} (\vec r) = \tau g_{\scalebox{0.6}{\ensuremath N}}  \, ,
\end{equation}
where the Lagrange multiplier $\tau$ represents the repulsive normal force between the two rods, and the non-positive gap function~$g_{\scalebox{0.6}{\ensuremath N}}$ depends on the minimum distance between the two rods. The distance between the lines of centroids of the two rods $\vec x^l (\theta^1_l)$ and $\vec x^ r (\theta^1_r)$ is given by 
\begin{equation}
d(\theta^1_l, \theta^1_r) = |\vec x^l (\theta^1_l) - \vec x^r (\theta^1_r)| \, , 
\end{equation}
and its minimum by
\begin{equation}\label{eq:minimumcontactd}
	d_{\min}= \mathop{\min}_{\theta^1_l, \theta^1_r} d(\theta^{1}_l, \theta^{1}_r) = d(\theta^{1c}_l, \theta^{1c}_r) \ . 
\end{equation}
This minimum can be determined using the Newton-Raphson scheme. Ultimately, the non-positive gap function~$g_{\scalebox{0.6}{\ensuremath N}}$ is defined as 
\begin{equation}
g_{\scalebox{0.6}{\ensuremath N}} = 
\begin{cases}
	d_{\text{min}}  - 2 r   & d_{\text{min}} < 2r \ , \\[5pt]
	0 & d_{\text{min}} \ge 2r   \ ,
\end{cases}
\end{equation}
where $r$ is the radius of the rods. 

The first variation of the contact potential energy~\eqref{eq:contact} gives its contribution to the weak form of the equilibrium equations~\eqref{eq:variation}. This contribution takes the form 
\begin{equation}\label{eq:deltacontact}
\delta \Pi_{c} = \tau \delta g_{\scalebox{0.6}{\ensuremath N}} + \delta \tau g_{\scalebox{0.6}{\ensuremath N}} \ .
\end{equation}
For further details on our contact implementation we refer to Herath~\cite{herathThesis2020}.

\subsection{Verification of the rod model \label{subsubsec:rodmodelVerif}} 
%
\begin{figure}[] 
\centering
\subfloat[]{
\includegraphics[scale=1]{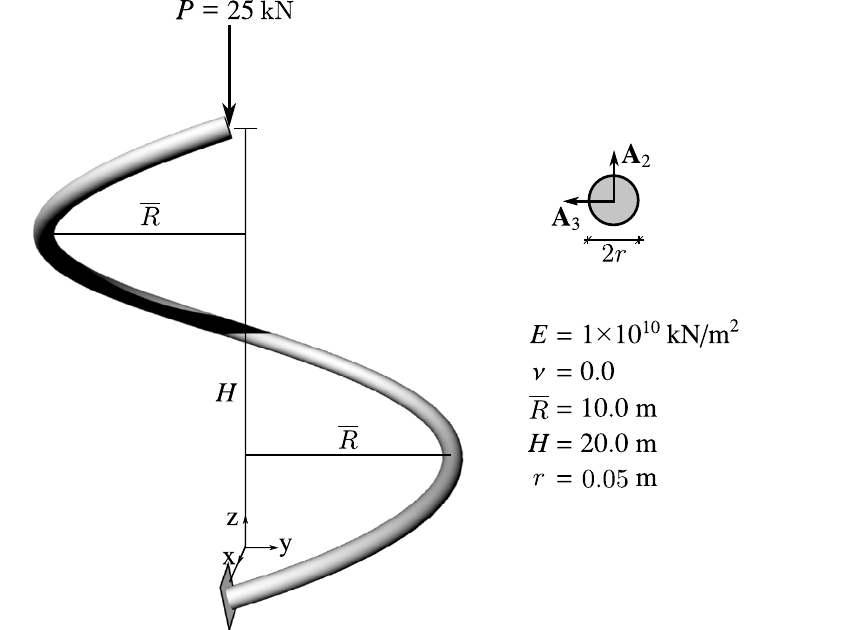}
\label{fig:rod_Spring_ProblemDescription}
}
\subfloat[]{
\advance\leftskip-4cm
\includegraphics[scale=0.5]{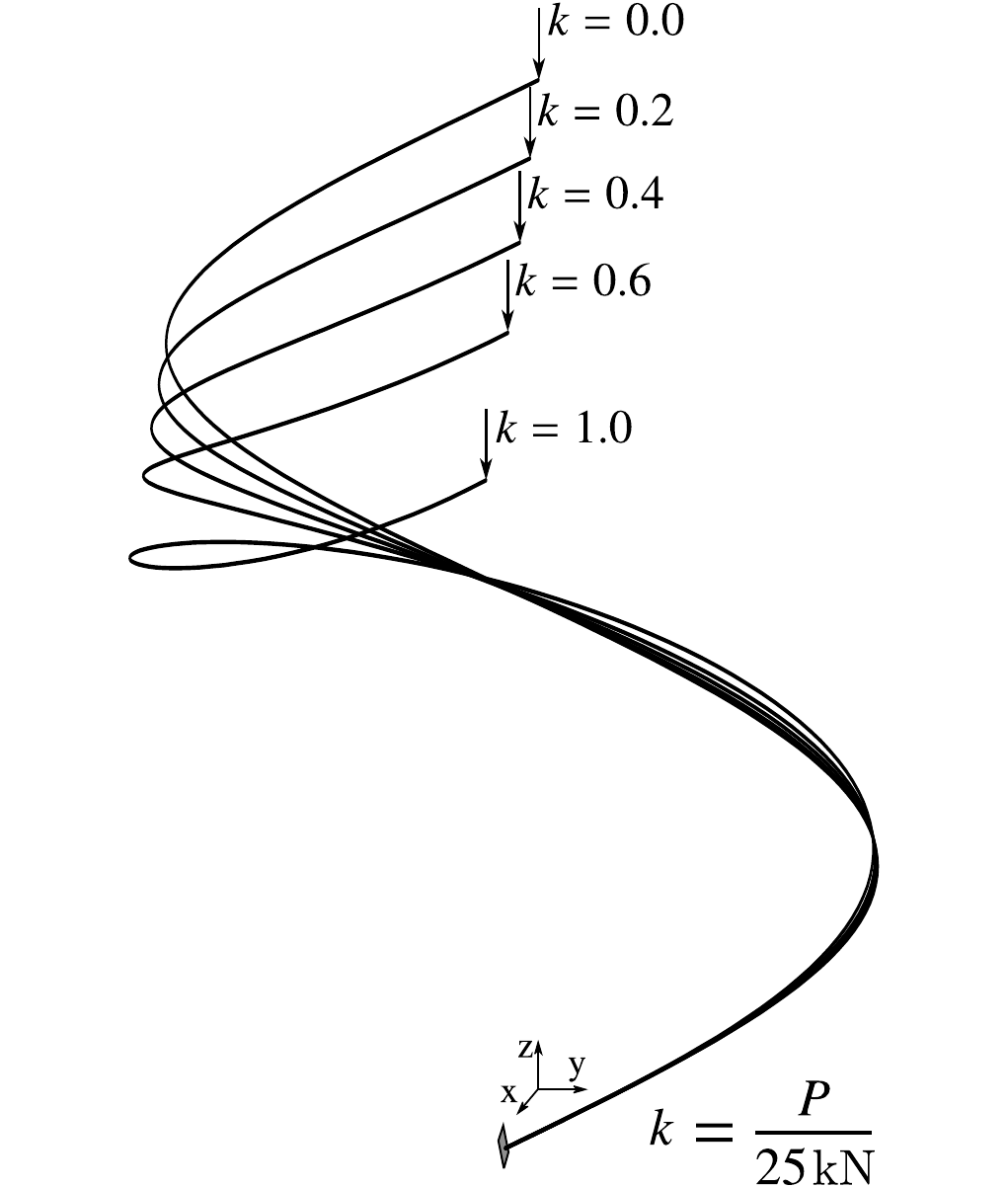}
\label{fig:rod_Spring_Deformed}
}
\caption{One-sided clamped helicoidal spring with tip loading. Problem description (a) and deformed shapes under different loads (b).}
\label{fig:rod_Spring}
\end{figure}
\begin{figure}[]
\centering
\includegraphics[scale=0.7]{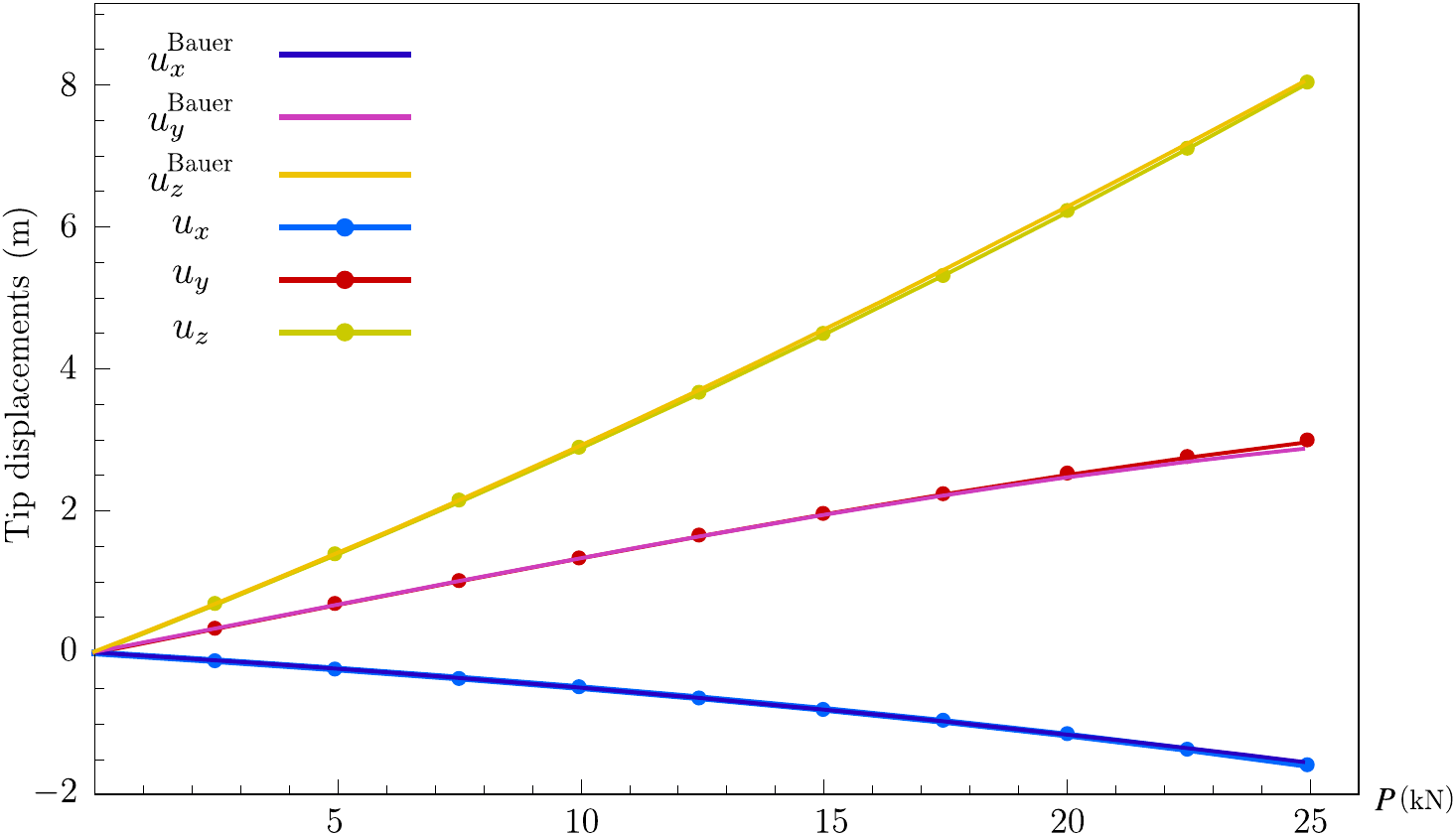}
\caption{Comparison of the obtained tip displacements of the helicoidal spring with Bauer et al.~\cite{bauer2016nonlinear}.}
\label{fig:rod_Spring_Comparison}
\end{figure}
We consider a membrane-bending-torsion interaction problem to verify the accuracy of the presented rod formulation. A spatial helicoidal spring is clamped at one end and a vertical load of magnitude 25 kN is applied at the other end, see Figure \ref{fig:rod_Spring_ProblemDescription}. 
The reference geometry of the line of the centroids is given by
\begin{align}
\bm X(\theta^1) &= \{ 10 \sin(2 \pi\,\theta^1), 10 \cos(2 \pi\,\theta^1), 20\,\theta^1 \}, \qquad \theta^1 \in [0,1] \, .
\end{align}
The material and geometric properties of the spring including the reference director orientations~$\vec A_\iota$ are given in~Figure \ref{fig:rod_Spring_ProblemDescription}.  The tip load is increased from zero to 25 kN in 10 uniform load steps. In Figure~\ref{fig:rod_Spring_Deformed} the deflected shapes for five different load levels are shown. As can be seen in  Figure~\ref{fig:rod_Spring_Comparison}  the obtained tip displacements are in excellent agreement with the results presented in Bauer~et~al.~\cite{bauer2016nonlinear}. In Bauer~et~al.~\cite{bauer2016nonlinear}, the reference directors $\vec A_i$ are first aligned with the principal axes of the rod cross-section and then mapped to the deformed $\vec a_i$ using the Euler-Rodrigues rotation formula. In this way, non-symmetric cross-sections with arbitrary rotations can be handled more efficiently.

\section{Computational Homogenisation \label{sec:compHom}}
\subsection{Microscale analysis \label{subsubsec:microscaleanalysis}} 
\begin{figure}[] 
\includegraphics[width=0.6\textwidth]{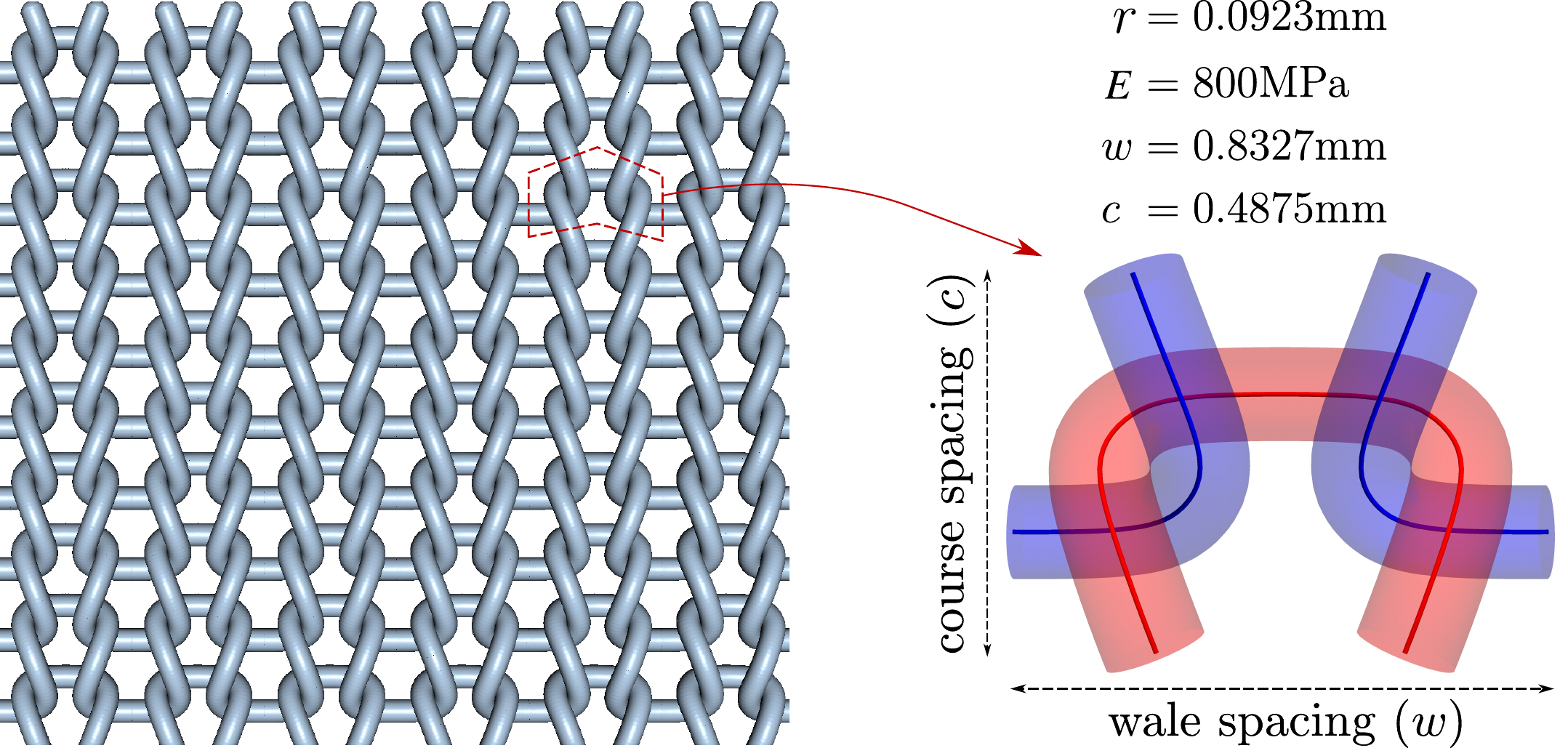}
\centering
\caption{Characteristic weft-knitted RVE selection and the definition of fibre geometric and material parameters.}
\label{fig:macroNmicro}
\end{figure}

%
A characteristic representative volume element (RVE) as depicted in~Figure \ref{fig:macroNmicro} is chosen to represent the periodic microstructure of a weft-knitted membrane. The intricate spatial arrangement of the yarns in the RVE is modelled, as in recent works~\cite{dinh2018,weeger2018}, using the approximate geometry proposed by Vassiliadis et al.~\cite{vassiliadis2007a}, see Appendix \ref{app:KnitGeometry}. For alternative geometric descriptions see~\cite{wadekar2020optimized} and references therein. First the yarn centrelines are defined and then the yarns are generated by sweeping a circle along those lines. 
On the macroscale we model the membrane with subdivision shell finite elements~\cite{Cirak:2001aa,Cirak:2011aa,long2012} and consider for homogenisation only the in-plane membrane response. The very small out-of-plane bending stiffness contribution of the rods to the membrane bending stiffness is not taken into account.

In first-order homogenisation, the macroscale membrane deformation gradient\ $\vec{F}_M$ is used to define the boundary conditions on the RVE. Here and in the following the macroscale and microscale quantities are denoted by subscripts $M$ and $m$, respectively. It is assumed that the volume averaged deformation gradient $\vec{F}_m$ of the rod in the RVE and the membrane deformation gradient~$\vec F_M$ are equal, that is, 
\begin{equation}
	\vec F_{M} = \dfrac{1}{V_{{\RVE}}} \int_{V_{{m}}} \vec F_{m} \D V \ ,
\label{eq:macromicroScaletransition}
\end{equation}
where $V_{\RVE}$ is the RVE volume and~$V_m$ is the rod volume within the RVE  in the reference configuration. Without going into details, the deformation gradient $\vec{F}_m = \partial \vec x / \partial \vec X$ of the line of centroids is given by the kinematic assumption~\eqref{eq:rodKinematics}. The Hill-Mandel lemma states that the macroscale and microscale work for an RVE must be equal
\begin{equation} 
	 \vec P_{M} \colon  \vec F_{M} =  \dfrac{1}{V_{{\RVE}}}  \int_{V_{m}}  \vec P_{m} \colon  \vec F_{m}  \D V  \, ,
\label{eq:PK1VolumAVGd}
\end{equation}
where~$\vec P_M$ is the macroscale first Piola-Kirchhoff membrane stress and~$\vec P_m$  is the microscale first Piola-Kirchhoff rod stress. The integral represents the internal virtual work of the rod within the RVE given by~\eqref{eq:DELTAinernalEnergy1} and~$\vec F_m$ is the respective virtual deformation gradient. For the lemma~\eqref{eq:PK1VolumAVGd} to hold only certain types of RVE boundary conditions can be chosen, including Dirichlet and periodic; see~\cite{geers2017,yvonnet2013computational} for details. Moreover, as a consequence of this lemma the macroscopic first Piola-Kirchhoff stress~$\vec P_M $ can be obtained from the volume averaged internal energy density  
\begin{subequations}
\begin{align}
	\vec P_M &= \frac{\partial \psi_M (\vec F_M) }{\partial \vec F_M}
	\intertext{with}
	 \psi_M (\vec F_M) &= \inf_{ \vec F_m \in \set K} \dfrac{1}{V_{{\RVE}}} \int_{V_{{m}}} \psi_m (\vec F_m) \D V \, , 
\end{align}
\end{subequations}
where the set~$\set K$ denotes the deformation gradients satisfying the RVE boundary conditions. The microscopic energy density~$\psi_m$ of a rod with an isotropic St Venant-Kirchhoff material is given by~\eqref{eq:stVenantEnergyDens}.

\begin{figure} 
\subfloat[]{
  \includegraphics[scale=0.65]{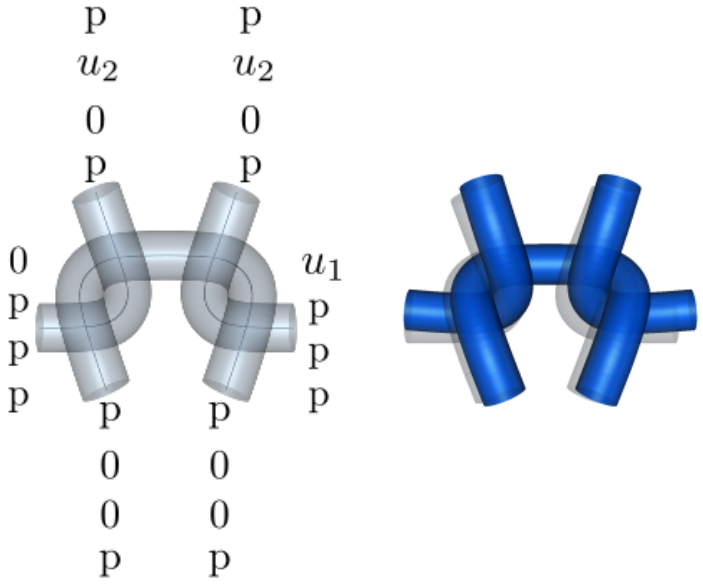}
}
\hspace{10mm}
\subfloat[]{
  \includegraphics[scale=0.65]{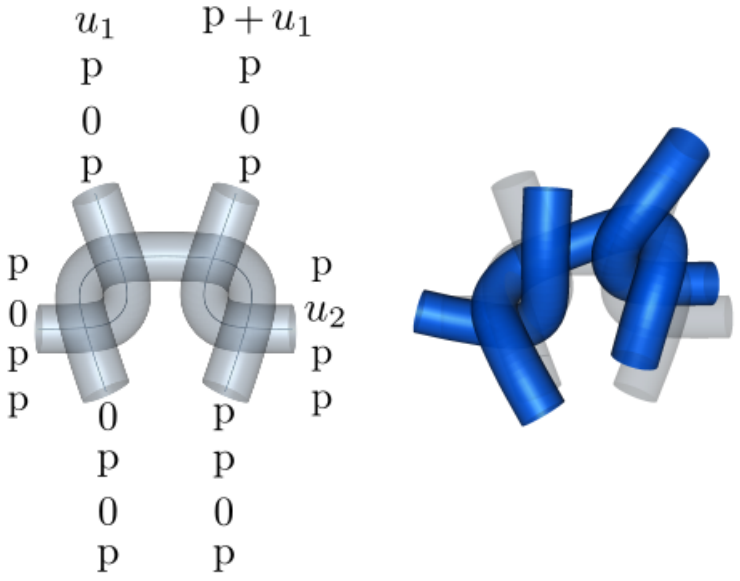}
}
\centering
\caption{Prescribed periodic boundary conditions and the respective deformed configuration of an RVE subjected to biaxial tension (a) and shear (b). At the six boundary finite element nodes, the four labels describe the boundary conditions for the three displacements and the one rotation, i.e. $(u_1, \, u_2, \, u_3, \, \vartheta )$, where $u_i$ refers to the applied displacements in the $\vec e_i$ direction. The labels $p$ and $0$ denote periodic and zero displacement or twist constraints, respectively.}
\label{fig:knitRVEPBC}
\end{figure}

In this work we choose as RVE boundary conditions the periodic boundary conditions depicted in Figure \ref{fig:knitRVEPBC}. Due to the  orthotropy of  the RVE response, the biaxial stretching and shear response are decoupled. Therefore, we consider the two cases separately and choose for each different boundary conditions. Boundary displacements in the thickness direction $\vec e_3$ are zero to simulate plane stress conditions.   
%
\subsection{Verification and validation of the microscale model \label{subsubsec:microscaleValidation}} 
\begin{figure} 
\includegraphics[width=0.7\textwidth]{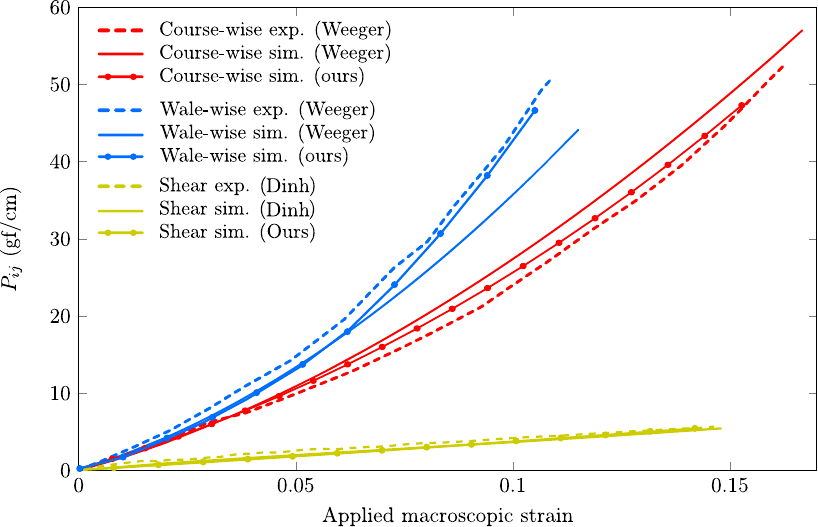}
\centering
\caption{Comparison of the obtained homogenised stresses for an RVE under course-wise and wale-wise uniaxial tension and pure shear with the experimental and numerical results reported in Weeger et al.~\cite{weeger2018} and Dinh et al.~\cite{dinh2018}. The unit conversion for stresses is $1.00\,\mathrm{N/mm} = 1019.716213\,\mathrm{gf/cm}$.}
\label{fig:microscaleVerification}
\end{figure}

%
To verify and validate our microscale model we consider RVEs under wale-wise and coarse-wise uniaxial tension and shear as reported in Dinh et al.~\cite{dinh2018} and Weeger et al.~\cite{weeger2018}. For a detailed problem description we refer to~\cite{dinh2018} and~\cite{weeger2018}. In our computations the yarn is discretised with~128 rod finite elements. In Dinh et al.~\cite{dinh2018} the yarn is modelled as a 3D solid and in Weeger et al.~\cite{weeger2018} as a discretised Euler-Bernoulli rod with the collocation method. As evident from Figure~\ref{fig:microscaleVerification} our results are in good agreement with the experiments and they appear to be be better than the numerical results by others. Throughout this paper, we use the yarn geometry and material parameters given in the above-mentioned two papers.
%

\section{Data-driven homogenisation\label{sec:knitGPRhomogen}}
\subsection{Review of Gaussian process regression}
Gaussian process regression is a statistical inference method rooted in Bayesian statistical learning. In the following, we briefly review the key steps in Gaussian process regression. For further details we refer to Rasmussen et al.~\cite{rasmussen2006}.
To begin with, we assume as a prior in the Bayesian formulation a random function $f(\vec{z})$ given by the zero-mean Gaussian process
\begin{equation}\label{eq:GPRdef}
f(\vec z) \thicksim \set G \set P \left( 0, k(\vec z, \vec z') \right) \, .
\end{equation}
The respective covariance function $k(\vec z, \vec z')$ is chosen to be the (stationary and isotropic) squared exponential function 
\begin{align}\label{eq:GaussianCov}
k(\vec z , \vec z') = \sigma^2_f \text{exp} \left( -\dfrac{(\vec z - \vec z')^\trans (\vec z - \vec z') }{2\ell^2} \right) \ ,
\end{align}
where $\sigma_f$ and $\ell$ are the yet to be determined scaling and characteristic lengthscale hyperparameters. For convenience, we define $\vec \Theta = (\sigma_f, \, \ell)$ as the vector of hyperparameters. The infinitely smooth covariance function~\eqref{eq:GaussianCov} encodes our prior assumptions about~$f(\vec z)$ before observing any training data. Other choices are possible, see \cite[Chapter~4]{rasmussen2006}.

Next, we consider a database comprised of a known training dataset and testing points. Each data point is given by an input vector $\vec{z}_i$ and a corresponding scalar observation $y_i$, with $i=1\,,2\,,\cdots\,, n_{GP}$. All training data points are collected in ($\mat Z$, $\mat y$) and all testing data points with unknown $\mat y_*$ in ($\mat Z_*$, $\mat y_*$). Moreover, we define the covariance matrix $\mat K$ with the components $\mathsf{K}_{ij} = k(\vec z_i, \vec z_j)$.
The discretisation of the Gaussian process~\eqref{eq:GPRdef} for the considered training and test data is given by the multivariate Gaussian distribution 
\begin{equation}\label{eq:GPpredictions}
\begin{bmatrix}
    \mat y \\
    \mat y_* 
\end{bmatrix} \sim \set N \left( \vec 0,
\begin{bmatrix}
\mat K(\mat Z,\mat Z) & \mat K(\mat Z,\mat Z_*) \\
 \ \mat K(\mat Z_*,\mat Z) & \ \mat K(\mat Z_*,\mat Z_*) 
\end{bmatrix}
    \right) \ .
\end{equation}
The probability density of the unknown data $\mat y_*$ at the given locations $\mat Z_*$ conditioned on the known training data $(\mat Z, \mat y)$ reads
\begin{equation}\label{eq:PredictionConditioned}
\mat y_* | \mat Z_*, \mat Z, \mat y \sim \set N \left (  \mat K(\mat Z_*,\mat Z) \mat K(\mat Z,\mat Z)^{-1} \mat y, \mat K(\mat Z_*,\mat Z_*) - \mat K(\mat Z_*,\mat Z) \mat K(\mat Z,\mat Z)^{-1} \mat K(\mat Z,\mat Z_*) \right ) \ .
\end{equation}
Hence, the best estimate for $\mat y_*$ is given by the expectation 
\begin{equation}
\label{eq:PredictionMean}
\overline{\mat y}_* = \mat K(\mat Z_*,\mat Z) \mat K(\mat Z,\mat Z)^{-1} \mat y \ ,
\end{equation}
and the uncertainty of the estimate $\overline{\mat K}_*$ is given by the covariance matrix
\begin{equation}
\label{eq:PredictionVariance}
 \text{cov}(\mat y_*) = \mat K(\mat Z_*,\mat Z_*) - \mat K(\mat Z_*,\mat Z) \mat K(\mat Z,\mat Z)^{-1} \mat K(\mat Z,\mat Z_*).
\end{equation}
The covariance matrix~$\mat K(\mat Z,\mat Z)$ of size~$n_{GP} \times n_{GP}$  is dense for the covariance kernel~\eqref{eq:GaussianCov}. If the number of training points~$n_{GP}$ become too large to invert~$\mat K(\mat Z,\mat Z)$, alternative approaches leading to sparse covariance matrices need to be considered~\cite[Chapter~8]{rasmussen2006}.
%

%
%
To estimate the hyperparameters $\vec \Theta$, we consider the marginal likelihood, or the evidence,
\begin{equation}\label{eq::marginalLikelihood}
p(\mat y | \mat Z) = \int p(\mat y | \mat Z, \mat{f}) p(\mat{f} | \mat Z) \D \mat{f} \, ,
\end{equation}
with the likelihood $p(\mat y | \mat Z, \, \mat f)$ and the prior $p(\mat f | \mat Z) = \set N (\vec 0, \, \mat K(\mat Z, \, \mat Z))$, where $\mat f$ are function values evaluated at $\mat Z$, c.f.~\eqref{eq:GPRdef}. The marginal likelihood is the probability of observing the data $\mat y$ for given hyperparameters $\vec \Theta$. Hence, it suggests itself to choose the hyperparameters so that the probability observing the given training data~$\mat y$ is maximised. In practice, it is numerically more stable to compute the maximum of the log marginal likelihood given by 
\begin{equation}\label{eq::marginallikelihoodSimpliFied}
	\log \ p(\mat y | \mat Z) = -\dfrac{1}{2} \mat y^\trans \mat K^{-1} \mat y - \dfrac{1}{2} \log |\mat K| - \dfrac{n_{GP}}{2} \log 2 \pi \ .
\end{equation}
Note that the $\log$ is a monotonic function so that 
\begin{equation} \label{eq:maxTheta}
	\vec \Theta^* =  \argmax_{\vec \Theta}  \ p(\mat y | \mat Z)   = \argmax_{\vec \Theta}  \left ( \log \ p(\mat y | \mat Z) \right )   \, .
\end{equation}
We use the scikit-learn Python library~\cite{scikitlearn} to find~$\vec \Theta^*$ with the gradient-based optimisation algorithm L-BFGS. The marginal likelihood is a non-convex function and certain care has to be taken to find the global maximum. After the optimal hyperparameters are determined they are used in computing the mean~\eqref{eq:PredictionMean} and covariance~\eqref{eq:PredictionVariance}.

\subsection{Gaussian process homogenisation \label{subsec:GPRcompHomog}}
%
\begin{figure} 
\includegraphics[width=0.8\textwidth]{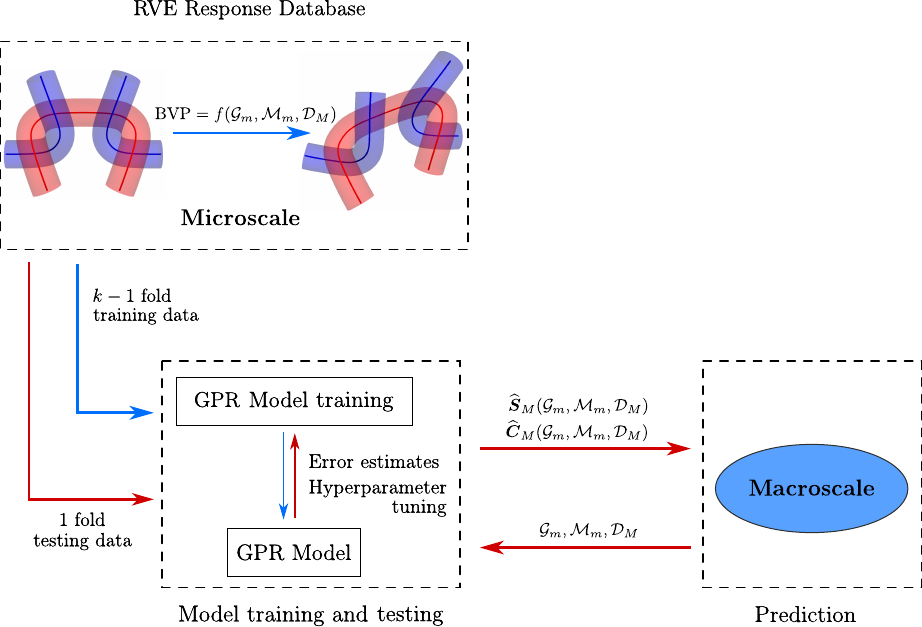}
\centering
\caption{Data-driven computational homogenisation and material design framework using Gaussian process regression.}
\label{fig:dataDrivenFramework}
\end{figure}
The implemented data-driven homogenisation framework is illustrated in Figure \ref{fig:dataDrivenFramework} and closely follows Bessa et al.~\cite{bessa2017}.
Data-driven homogenisation begins with constructing a response database for the RVE. This response database is designed by defining the microscale boundary value problem (BVP) using three types of design variables, namely geometric properties $\set G_m$, material properties $\set M_m$ and boundary conditions $\set D_M$. In this paper, we concentrate on the design variable $\set D_M$ relevant to homogenisation. It is straightforward to include $\set G_m$ and $\set M_m$ which lead to the material design of knitted membranes. 

Next, we  determine the hyperparameters of the GPR model using $k$-fold cross-validation rather than directly using the~$\vec \Theta^*$ obtained with~\eqref{eq:maxTheta} considering the entire training dataset. The GPR model is iteratively evaluated  by using $k-1$ folds for training and the remaining one for testing. We choose $k=5$ which provides a good compromise between accuracy and efficiency. As error metrics for the predicted strain energy and stress resultants we use the correlation of determination $\mathcal{R}^2$ and the mean squared error (MSE) given by
\begin{subequations}
\begin{align}
\mathcal{R}^2(\mat y,\widehat{\mat y}) &= 1 - \dfrac{\sum^{n_T}_{i=1}(y_i - \widehat{y}_i)^2}{\sum^{n_T}_{i=1}(y_i - \overline{y})^2}\label{eq:R^2} \, ,\\
\text{MSE}(\mat y,\widehat{\mat y}) &= \dfrac{1}{n_T} \sum^{n_T}_{i=1} (y_i - \widehat{y}_i)^2 \label{eq::MSE} \, ,
\end{align}
\end{subequations}
where $n_T = \lfloor{n_{GP} / (k-1)} \rfloor$ is the number of the test points, $\mat y$ is the target output, $\widehat{\mat y}$ is the predicted output and $\overline{y}$ is the expectation of the target output given by the GPR model. When evaluating the MSE of the predicted stress tensor, we take the square root of the squared sum of the MSE of each stress resultant
\begin{align}\label{eq::MSEofS}
\text{MSE}(\vec S, \widehat{\vec{S}}) = \sqrt{\text{MSE}^2(S_{11},\widehat{S}_{11}) + \text{MSE}^2(S_{22},\widehat{S}_{22}) +\text{MSE}^2(S_{12},\widehat{S}_{12})}.
\end{align}
The GPR model training is implemented in a Python environment using scikit-learn machine learning library~\cite{scikitlearn}. For each test fold, we determine the corresponding hyperparameters by maximising the log marginal likelihood~\eqref{eq:maxTheta}. It is expected that the obtained hyperparameters have similar values for all test folds. In practice, this is not the case because the log marginal likelihood is usually a non-convex function.  Therefore, we vary during optimisation the upper and lower limits of the hyperparameters so that all test folds yield similar hyperparameter values.  Once this is achieved, we record the hyperparameters~$\vec \Theta^*$ of the GPR model, which has the lowest MSE, and save the model for subsequent macroscale analyses.

Lastly, we integrate our in-house thin-shell solver~\cite{cirak2000, long2012} with the trained GPR model to simulate the homogenised knitted membrane. As shown in Figure \ref{fig:dataDrivenFramework}, for given macroscale deformations (strains) the trained GPR model yields the corresponding predicted stress resultants and their tangents to the macroscale shell solver.

\subsection{Macroscale analysis \label{subsec:macroscaleAnalysis}} 
%
For plane-stress membrane deformation, a response database is created with the three in-plane components of the Green-Lagrange strain~$\vec E$ as the inputs and RVE volume averaged strain energy $\psi (\vec E)$ as the output, i.e. $\vec z_i$ and~$y_i$ in Section~\ref{subsec:GPRcompHomog}. With a slight abuse of notation~$\vec E$ is here, in contrast to Section~\ref{subsec:rodkinematics}, a two-dimensional second order tensor. In the four-dimensional response database, points have the coordinates \mbox{$(\vec E , \psi) = (E_{11},E_{22},E_{12}, \psi)$.} The strain components $E_{11}$ and $E_{22}$ correspond to a biaxial deformation state that is decoupled from the shear deformation state with $E_{12}$, see Section~\ref{subsubsec:microscaleanalysis}. The RVE strain energy~$\psi$ is obtained by considering a biaxial strain state with~$E_{11}$ and $E_{22}$ and a shear deformation state with~$E_{12}$ and summing up their strain energies. We use the validated strain limits of Figure~\ref{fig:microscaleVerification} in constructing the response database. Thus, allowing for a 5$\%$ compressive strains, we define the design variables of the response database as,
\begin{align}\label{eq:dataDrive_D_M}
\set D_{M} &=
		 \begin{cases}
            \set D^1_{M} = E_{11} \in [-0.05,0.15] \, , \\
            \set D^2_{M} = E_{22} \in [-0.05,0.15] \, , \\
            \set D^3_{M} = E_{12} \in [0.00,0.15] \, .
        \end{cases}
\end{align}
We solve the RVE problem with the boundary conditions stipulated by the strain states defined in \eqref{eq:dataDrive_D_M} and store the strain components and respective volume averaged energies in the response database.  After passing the GPR model training and testing phase, as shown in Figure \ref{fig:dataDrivenFramework}, the strain energy density~$\widehat \psi$ of the plane-stress deformation is predicted for a given deformed state of the RVE. Plane stresses and stress tangents are computed by differentiating~\eqref{eq:PredictionMean}, that is,
\begin{equation}
\widehat{\vec{S}}_M = \dfrac{\partial \widehat{\psi}_M}{\partial \vec E_M} \qquad \text{and} \qquad \widehat{\vec{C}}_M = \dfrac{\partial \widehat{\vec{S}}_M}{\partial \vec E_M} \, .
\end{equation}

\section{Examples \label{sec:examples}}
\subsection{GPR model training and testing \label{subsec:GPRresults}}
We start with presenting the error metrics of GPR model training and testing for varied sizes of response databases. GPR model training errors are observed to take values \mbox{$\mathcal{R}^2 = 1.0$} and \mbox{$\text{MSE } < 10^{-8}$.} Hence, we present only the error metrics of testing datasets. All the following computations are performed on the Intel Core i5-4590 CPU @ 3.30GHz $\times$ 4 processor.

\begin{figure} 
\includegraphics[scale=0.65]{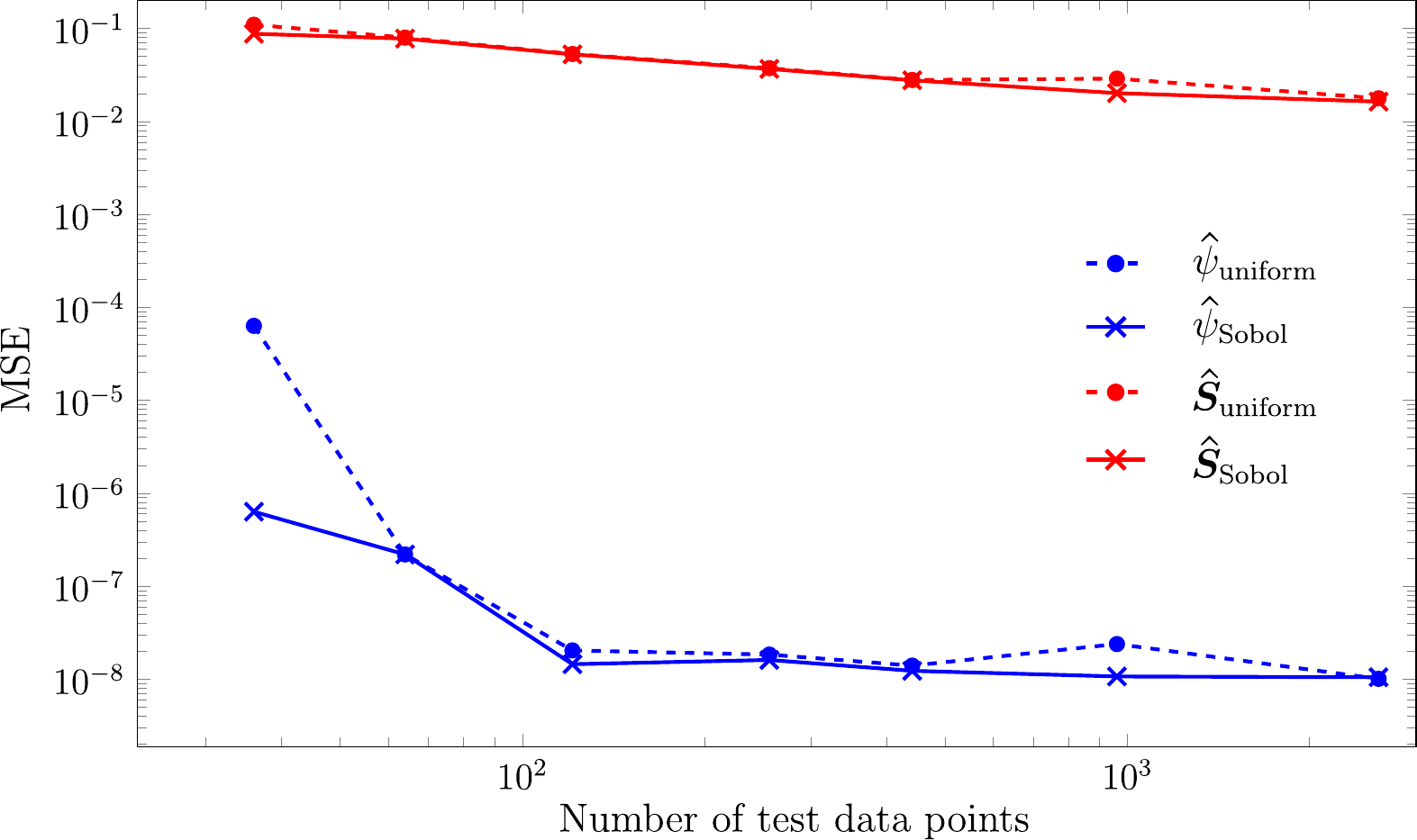}
\centering
\caption{GPR model training and testing. The mean squared error of the predicted strain energy density and the second Piola-Kirchhoff stress resultants for uniform and Sobol sampling.}
\label{fig:GPRerrors}
\end{figure}

\begin{figure} 
\centering
\subfloat[]{
\includegraphics[width=0.45\columnwidth]{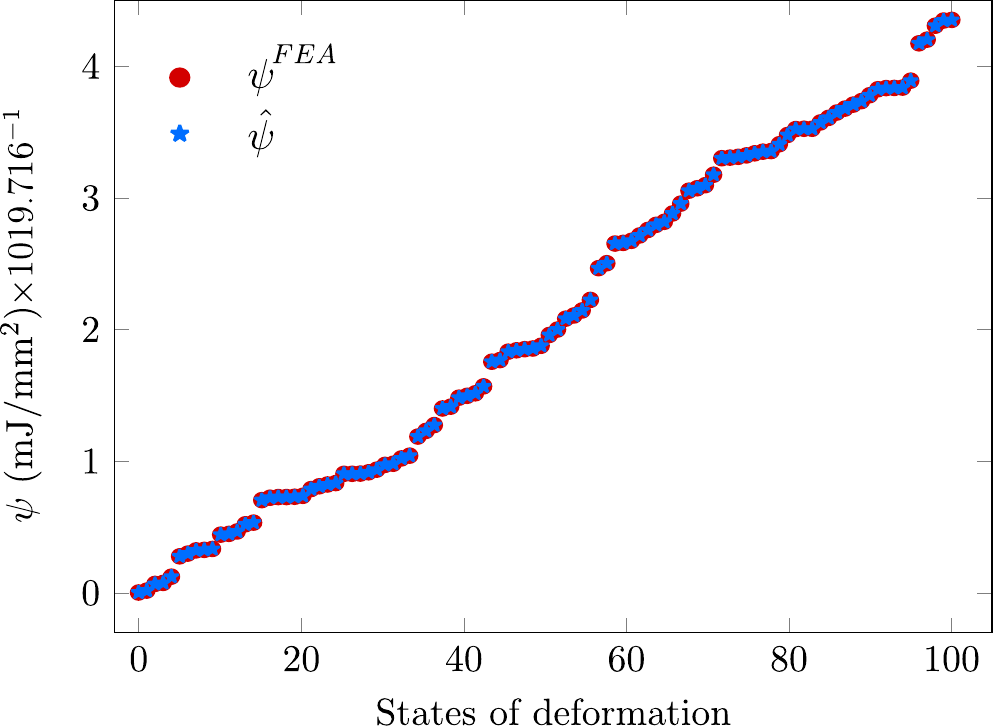}
\label{fig:energyPrediction12}
}
\hspace{8mm}
\subfloat[]{
\centering
\includegraphics[width=0.45\columnwidth]{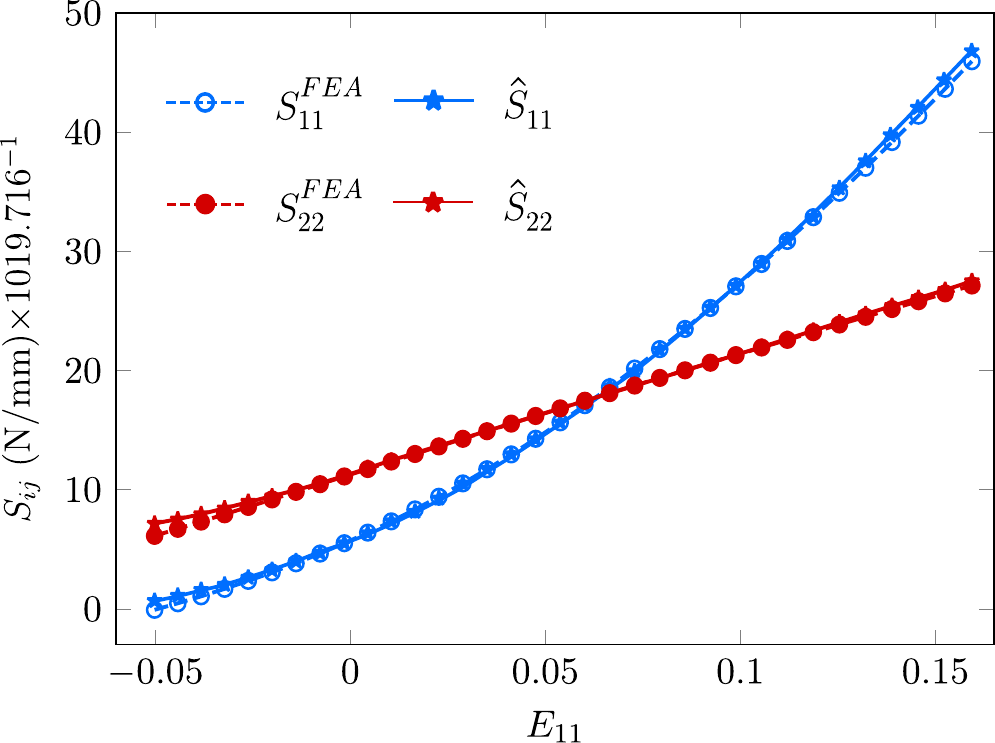}
\label{fig:stressPrediction12}
}
\centering
\caption{GPR model training and testing. Predictions of (a) strain energy density (sorted by ascending potential) and (b) stress resultants for a biaxial strain state of $E_{11} = [-0.05, 0.15]$, $E_{22} = 0.08$ and $E_{12} = 0$.}
\label{fig:stressPredictions}
\end{figure}

Figure \ref{fig:GPRerrors} shows the MSEs of predicted energy density and second Piola-Kirchhoff stress resultants, as well as a comparison between two sampling techniques, namely uniform sampling and Sobol sampling~\cite{sobol1976uniformly}, on the testing errors. Uniformly sampled data points are generated in a way that for the three features in the input vector $\vec z_i$, $n^3_S$ entries occupy the response database, where $n_S$ is the number of uniformly distributed entries in each feature. In total seven datasets with \mbox{$n_S^3 \in \{ 216, \, 343, \, 729, \, 1331, \, 1728, \, 2197, \, 4913, \, 12167\}$} are considered; $80 \%$ of each are taken as training points and the remaining as testing points. Moreover, $\mathcal{R}^2=1.0$ is rounded to the fifth decimal number for all response databases. MSEs in stress predictions are comparatively higher than those in energy predictions but remain less than 2$\%$. Considering the error convergence and algorithmic efficiency, we use a response database with 2601 test data points in the subsequent macroscale simulations. Furthermore, this particular trained GPR model has the optimum hyperparameter values $\vec{\Theta}^* = (\sigma_f, \ell) = (1.25213, 0.00904)$ and the maximum log marginal likelihood $ 9594.98391$. Model training time was recorded as 29 minutes and 48 seconds.
In Figure~\ref{fig:stressPredictions}, we visualise and compare the predictions of a subsample of the chosen response database. Figure \ref{fig:energyPrediction12} depicts the strain energy density predictions, whereas Figure \ref{fig:stressPrediction12} presents the stress predictions for the biaxial homogenised response of an RVE.
%
\subsection{Stretched membrane I: comparison of yarn-level and homogenised displacements  \label{subsec:yarnlevelPull}}
%
\begin{figure} 
\centering
\subfloat[]{
\includegraphics[width=0.48\columnwidth]{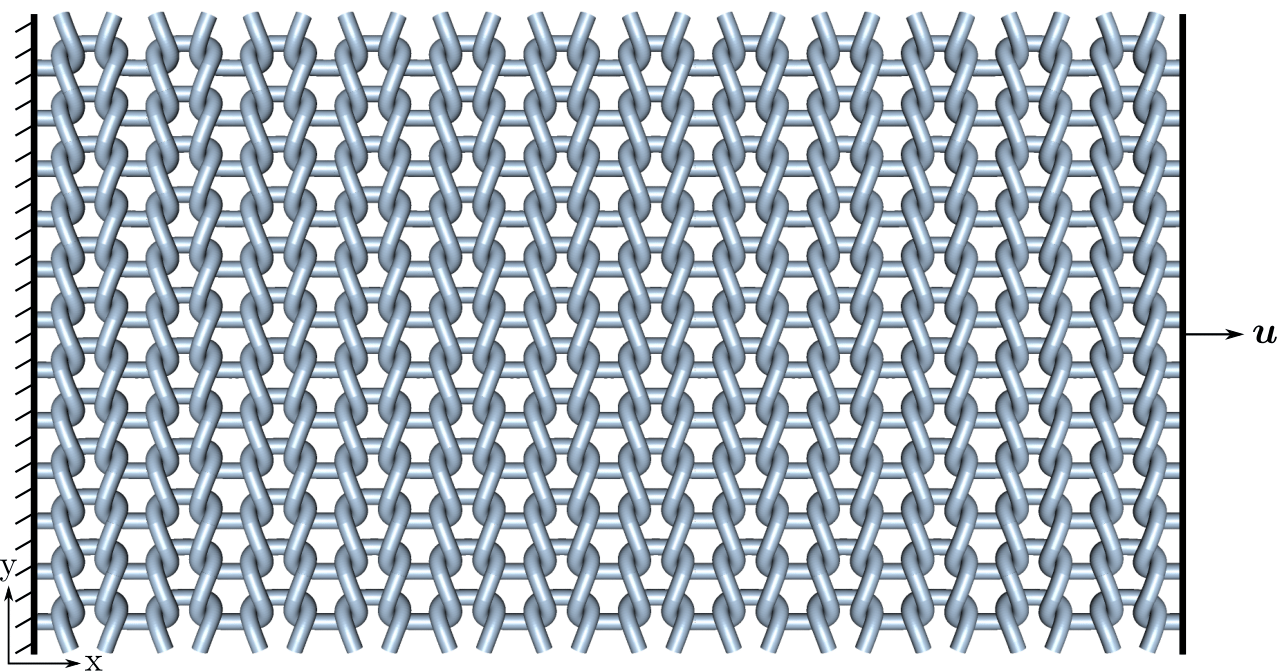}
\label{fig:yarnlevel_ProbDesc}
}
\centering
\subfloat[]{
\centering
\includegraphics[width=0.48\columnwidth]{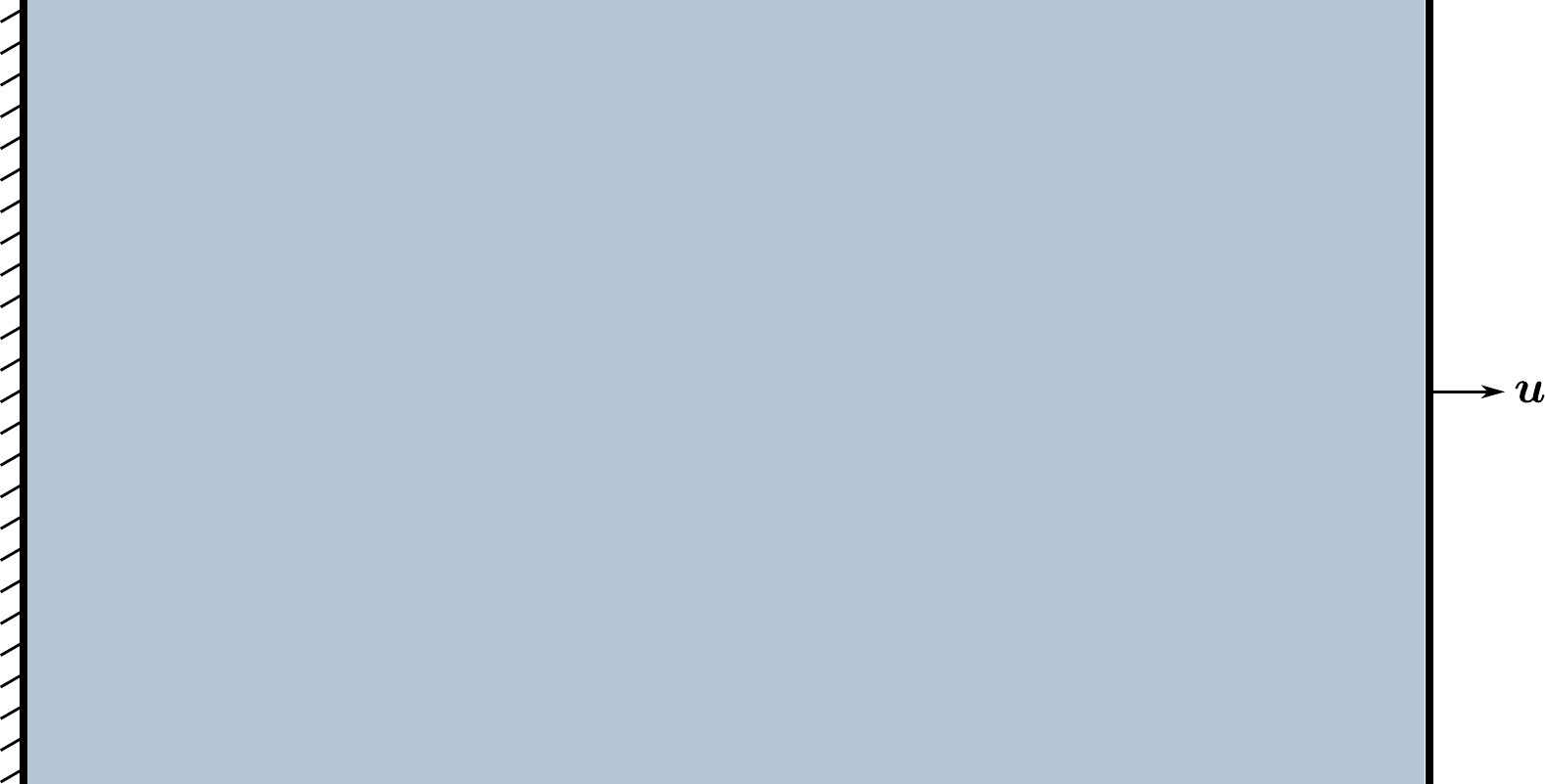}
\label{fig:homog_ProbDesc}
}
\centering
\caption{Stretched membrane I. Problem description for the yarn-level (a) and the homogenised membrane (b) models.}
\label{fig:KnitprobDesc}
\end{figure}
\begin{figure} 
\centering
\subfloat[]{
\includegraphics[width=0.49\columnwidth]{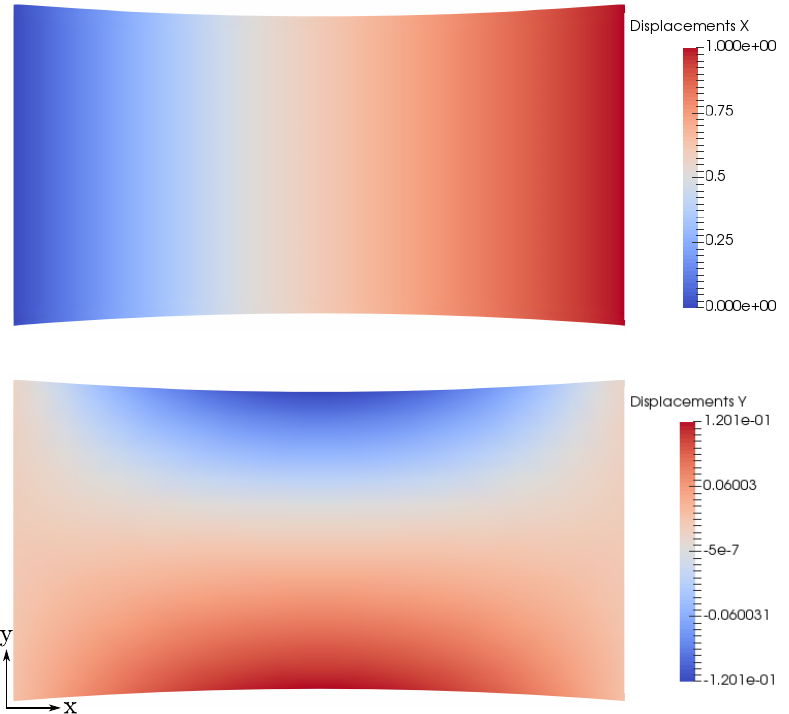}
\label{fig:SimulationyarnLevel}
}
\centering
\subfloat[]{
\centering
\includegraphics[width=0.49\columnwidth]{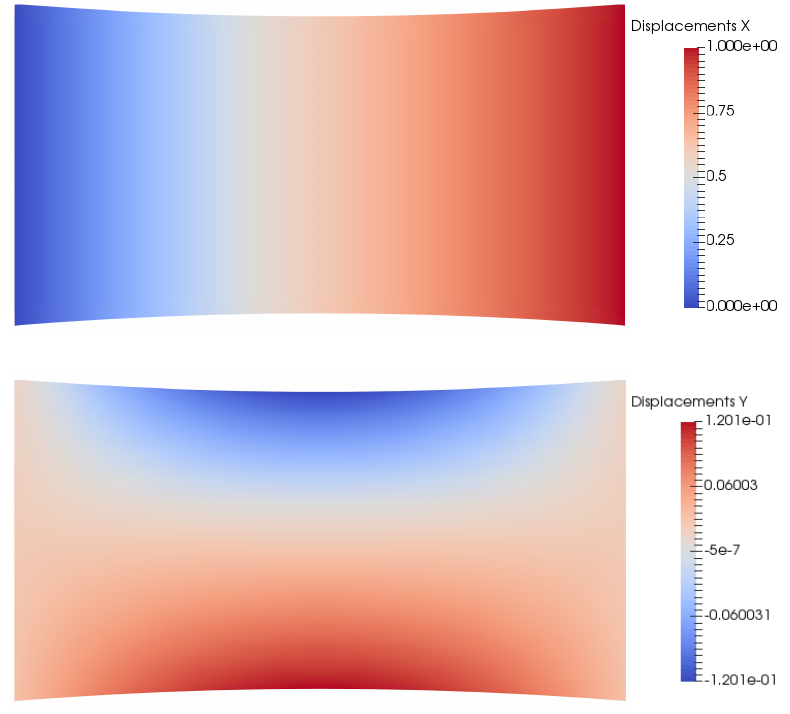}
\label{fig:SimulationHomogenised}
}
\centering
\caption{Stretched membrane I. Comparison of horizontal and vertical displacement iso-contours of the yarn-level (a) and homogenised membrane (b) models.}
\label{fig:yarnvsHomog}
\end{figure}
We consider a membrane under uniaxial tension and compare its response using the Gaussian process homogenised material model and an equivalent yarn-level model. The yarn-level model comprises 12 loops in the course direction and 12 loops in the wale direction. We use for the analysis with the homogenised material model a membrane of equal size, that is, $10.00\,\mathrm{mm}$ long, $5.85\,\mathrm{mm}$ wide and $0.2356\,\mathrm{mm}$ thick. The membrane is discretised with a structured quadrilateral mesh with  $12\times12$ elements. Problem descriptions of the two models are presented in Figure \ref{fig:KnitprobDesc}.

Both models are stretched by $1\,\mathrm{mm}$ in the $x$-direction and the resulting $x$ and $y$ displacements are presented in Figure \ref{fig:yarnvsHomog}. For comparison purposes, the results of the yarn level model in Figure~\ref{fig:SimulationyarnLevel} are visualised by projecting the nodal values of the yarn onto the $x-y$ plane and interpolating them on a Delaunay mesh. The $y$-direction displacement distribution of the yarn model in Figure \ref{fig:SimulationyarnLevel} is slightly different from that of the homogenised membrane in Figure \ref{fig:SimulationHomogenised}. This difference has two causes.  First, zero out-of-plane displacement boundary conditions on the top and bottom edges of the yarn-level model directly contribute to the observed difference. These boundary conditions are applied to simulate selvedge stitches that restrict the top and bottom yarns from being freely straightened during deformation. Secondly, the geometric asymmetry of the yarn-level model about the mid-horizontal plane has an influence on the overall response due to the relatively small number of loops used in both course and wale directions. Furthermore, a very similar Poisson effect is observed in both the yarn-level and homogenised membrane models. The absolute maximum \mbox{$y$-displacement} is recorded with $0.12\,\mathrm{mm}$ for a stretch by $1\,\mathrm{mm}$ in the $x$-direction.

\subsection{Stretched membrane II: homogenised stresses  \label{subsec:pulltests}}
\begin{figure} 
\centering
\subfloat[]{
\includegraphics[width=0.99\textwidth]{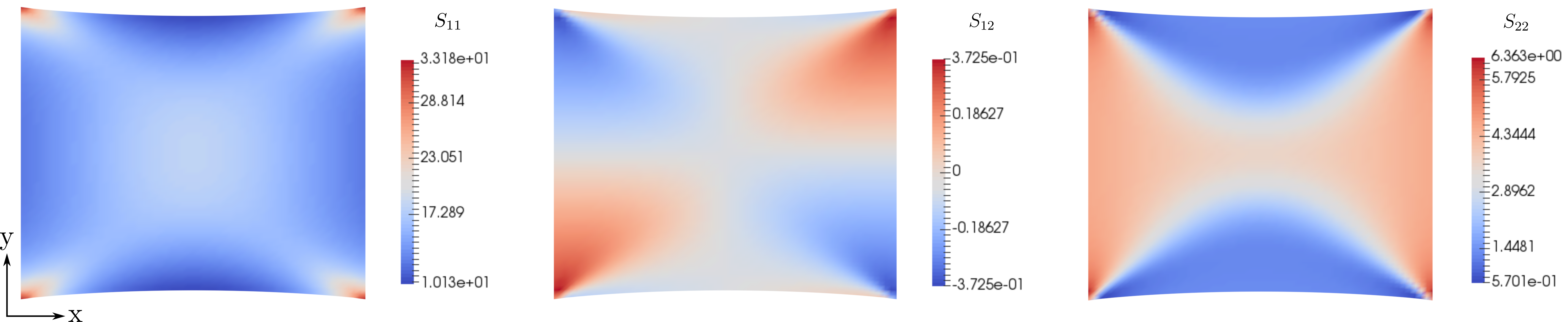}
\label{fig:pulltTestCoursewise}
}

\centering
\subfloat[]{
\includegraphics[width=0.99\textwidth]{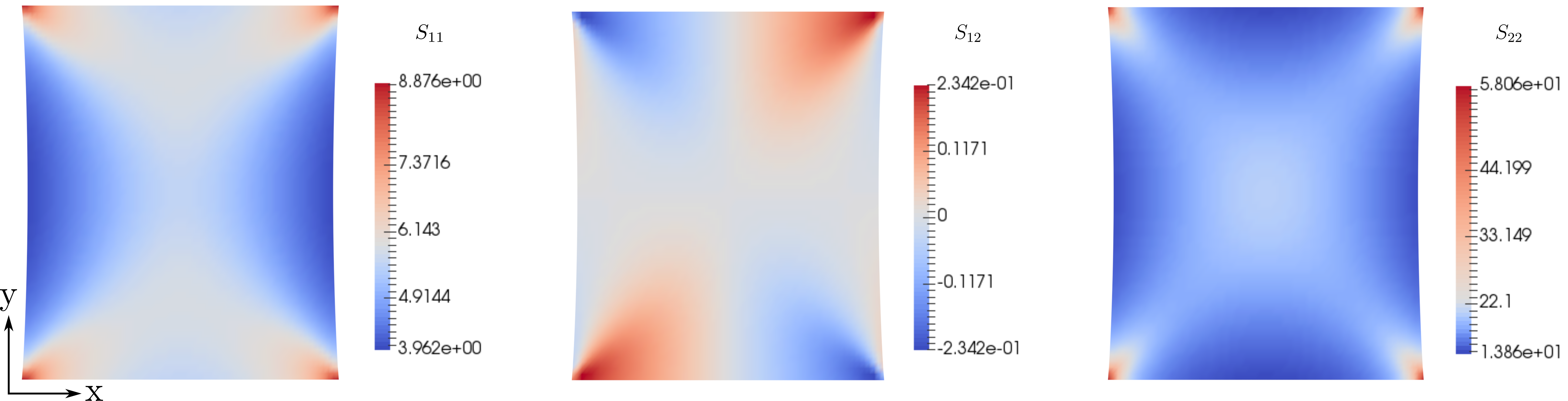}
\label{fig:pulltTestWalewise}
}
\centering
\caption{Stretched membrane II:  Second Piola-Kirchhoff stress contours of the homogenised  course-wise stretched (a) and wale-wise stretched (b) membranes.}
\label{fig:pulltTests}
\end{figure}

We perform two stretch tests on a square knitted membrane sheet with a side length of $10\,\mathrm{mm}$ in the course and wale directions. A structured quadrilateral mesh of size  $10 \times 10$ ($100$ elements) is used. During the displacement controlled deformation one edge is fixed while the other is stretched in the course or wale direction, respectively, by $1\,\mathrm{mm}$. The deformed shapes and stress contours are depicted in Figure \ref{fig:pulltTestCoursewise} for course-wise and in Figure \ref{fig:pulltTestWalewise} for wale-wise stretching. Figure \ref{fig:pulltTests} clearly shows the orthotropic response of the knitted membrane as the stress resultants are different depending on the direction of the stretching, comparing, e.g.,~$S_{11}$ in Figure~\ref{fig:pulltTestCoursewise} with~$S_{22}$ in Figure~\ref{fig:pulltTestWalewise}. Furthermore, the stiffer response in the wale-wise direction manifests itself in a higher stress $S_{22}$ in Figure \ref{fig:pulltTestWalewise} in comparison to $S_{11}$ in Figure \ref{fig:pulltTestCoursewise}.

\section{Conclusions \label{sec:conclusion}}
%
We introduced a data-driven approach for computational homogenisation of knitted membranes to address the challenges posed by conventional homogenisation schemes. The inherent large deformations in knitted textiles are accurately captured by the finite deformation rod model and compared to experimental and numerical results. The predictions of the microscale model under uniaxial tension and shear are in very good agreement with experimental results. As shown for uniaxially stretched membranes the homogenised model yields similar macroscopic displacements like a more elaborate full yarn-level model. Incorporating the statistical GPR model in computational homogenisation circumvents the need for expensive microscale simulations at each quadrature point of the macroscale membrane, thus yielding significant computational savings.

The presented approach can be further validated and extended in several ways. We considered in this paper only membranes subjected to in-plane tension with no out-of-plane deformations. Further experimental results are needed to validate the proposed data-driven approach in case of significant compression and bending, especially fine-scale wrinkling. Furthermore, we considered so far only weft-knitted membranes with a uniform stitch pattern. To produce complex three-dimensional surfaces, it is necessary to alter the stitch pattern, for instance, by locally increasing or decreasing the number of loops from row to row or  reducing the number of rows. Such changes represent discontinuities in the stitch pattern, and  suitable RVEs have to be defined to capture their homogenised response. Putting aside questions of RVE size and validity of homogenisation assumptions, it is straightforward to consider in the introduced data-driven approach RVEs with other stitch patterns. Furthermore, depending on the loading conditions knitted membranes are prone to geometric instabilities in form of wrinkling on the macroscale and rod buckling on the microscale. In presence of such instabilities the homogenisation assumptions are usually not valid and alternative approaches must be used~\cite{nezamabadi2009multilevel,miehe2002computational}. The data-driven variants of these approaches are essential for the analysis of large-scale knitted membranes. Finally, a key advantage of the data-driven approach is the prospect to consider design parameters, pertaining to the stitch geometry or yarn material, in the GPR model. This opens up the possibility to optimise those parameters with an efficient gradient-based optimisation algorithm. 

\appendix
\section{Derivatives of strains and internal forces}\label{app:intForces}
In this Appendix we summarise the detailed equations used in the implementation of the introduced finite deformation rod finite element.
\subsection{Strain derivatives}\label{App:strainDer}

The derivative of the components~$E_{ij}$ of the strain tensor~\eqref{eq:allStrainComponents} with respect to the nodal displacements $\vec u_I$ are given by
\begingroup
\allowdisplaybreaks
\begin{subequations}
\begin{flalign}
\frac{\partial \alpha}{\partial \vec u_I} & = \vec{a}_1 B^I_{,1} \label{eq:MembranestrainGradient} \, , &\\
\frac{\partial \beta_{\iota}}{\partial  \vec u_I} & = -\frac{\partial \vec{a}_{\iota}}{\partial  \vec u_I} \cdot \vec{a}_{1,1} - \vec{a}_{\iota} \cdot \frac{\partial \vec{a}_{1,1}}{\partial \vec u_I} \label{eq:BendingstrainGradient} \, , &\\
\frac{\partial \gamma}{\partial  \vec u_I} & = \frac{1}{2} \left( \frac{\partial \vec{a}_2}{\partial  \vec u_I} \cdot \vec{a}_{3,1} + \vec{a}_2 \cdot \frac{\partial \vec{a}_{3,1}}{\partial \vec u_I} \right) \, , \label{eq:TorsionstrainGradient}
\end{flalign}
\end{subequations}
\endgroup
where
\begingroup
\allowdisplaybreaks
\begin{flalign}
\frac{\partial \vec{a}_1}{\partial \vec u_I} &= B^I_{,1} \vec{I} \, , \nonumber&\\
\frac{\partial \vec{a}_{1,1}}{\partial \vec u_I} &= B^I_{,11} \vec{I} \, , \nonumber&\\
\frac{\partial \vec{a}_{\iota}}{\partial \vec u_I} &=
\frac{\partial \vec{\Lambda}_2}{\partial \vec u_I}\vec{\Lambda}_1{\vec A}_{\iota} \, , \nonumber&\\
\vec a_{\iota,1} &= \vec{\Lambda}_{2,\text{1}}\vec{\Lambda}_1\vec A_{\iota} +  \vec{\Lambda}_2\vec{\Lambda}_{1,\text{1}}\vec A_{\iota} +  \vec{\Lambda}_2\vec{\Lambda}_1\vec A_{\iota,1} \, , \nonumber&\\
\frac{\partial \vec{a}_{{\iota},1}}{\partial \vec u_I} &=
\frac{\partial \vec{\Lambda}_{2,\text{1}}}{\partial \vec u_I}\vec{\Lambda}_1{\vec A}_{\iota} +
\frac{\partial \vec{\Lambda}_2}{\partial \vec u_I}\vec{\Lambda}_{1,\text{1}}{\vec A}_{\iota} + 
	\frac{\partial \vec{\Lambda}_2}{\partial \vec u_I}\vec{\Lambda}_1{\vec A}_{{\iota},1} \, . \nonumber&
\end{flalign}
\endgroup
The derivative of the rotation matrix $\vec{\Lambda}_2$ in \eqref{eq:R2} with respect to the nodal displacements $\vec u_I$ reads
\begingroup
\allowdisplaybreaks
\begin{flalign}
\frac{\partial \vec{\Lambda}_2 }{\partial \vec u_I} &= \frac{\partial \vec{\Lambda}_2}{\partial \hat{\vec a}_1} \frac{\partial \hat{\vec a}_1}{\partial \vec u_I} \, , &
\end{flalign}
\endgroup
where
\begingroup
\allowdisplaybreaks
\begin{flalign}
\frac{\partial \vec{\Lambda}_2}{\partial \hat{\vec a}_1} &= -\frac{(\hat{\vec A}_1 + \hat{\vec a}_1) \otimes \vec{I}}{(1+ \hat{\vec a}_1 \cdot \hat{\vec A}_1)} - \frac{\vec{I} \otimes \hat{\vec a}_1}{(1+ \hat{\vec a}_1 \cdot \hat{\vec A}_1)} + \frac{((\hat{\vec A}_1 + \hat{\vec a}_1) \otimes \hat{\vec A}_1) \otimes \hat{\vec a}_1 }{(1+ \hat{\vec a}_1 \cdot \hat{\vec A}_1)^2} \, , \nonumber&\\
\frac{\partial \hat{\vec a}_1}{\partial \vec u_I} &= \left( \frac{\vec{I}}{|\vec a_{1}|} - \frac{\vec a_1 \otimes \vec a_1}{|\vec a_{1}|^3} \right) B_{,1}^I \, . \nonumber&
\end{flalign}
\endgroup
The derivatives of the bending and torsional shear strains in~\eqref{eq:allStrainComponents} with respect to the nodal twist $\vartheta_I$ are given by
\begingroup
\allowdisplaybreaks
\begin{subequations}
\begin{flalign}
\frac{\partial \beta_{\iota}}{\partial  \vartheta_I} & = -\frac{\partial \vec{a}_{\iota}}{\partial  \vartheta_I} \cdot \vec{a}_{1,1} \, , &\\
\frac{\partial \gamma}{\partial  \vartheta_I} & = \frac{1}{2} \left( \frac{\partial \vec{a}_2}{\partial  \vartheta_I} \cdot \vec{a}_{3,1} + \vec{a}_2 \cdot \frac{\partial \vec{a}_{3,1}}{\partial \vartheta_I} \right) \, , &
\end{flalign}
\end{subequations}
\endgroup
where
\begingroup
\allowdisplaybreaks
\begin{flalign}
\frac{\partial \vec{a}_{\iota}}{\partial \vartheta_I} &=
\vec{\Lambda}_2 \frac{\partial \vec{\Lambda}_1}{\partial \vartheta_I}{\vec A}_{\iota} \, , \nonumber&\\
\frac{\partial \vec{a}_{{\iota},1}}{\partial \vartheta_I} &=
\vec{\Lambda}_{2,\text{1}}\frac{\partial \vec{\Lambda}_1}{\partial \vartheta_I}{\vec A}_{\iota}
 + \vec{\Lambda}_2 \frac{\partial \vec{\Lambda}_{1,\text{1}}}{\partial \vartheta_I}{\vec A}_{\iota} + 
\vec{\Lambda}_2\frac{\partial \vec{\Lambda}_1}{\partial \vartheta_I}{\vec A}_{{\iota},1} \, . \nonumber&
\end{flalign}
\endgroup
The derivative of the rotation matrix $\vec{\Lambda}_1$ in \eqref{eq:R1}  with respect to the nodal twist $\vartheta_I$ is given by
\begingroup
\allowdisplaybreaks
\begin{subequations}
\begin{flalign}
\frac{\partial \vec{\Lambda}_1 }{\partial \vartheta_I} &=\frac{\partial \vec{\Lambda}_1 }{\partial \vartheta}\frac{\partial \vartheta}{\partial \vartheta_I} \, , \\
\frac{\partial \vec{\Lambda}_1 }{\partial \vartheta} &=  \cos\vartheta \vec{L} + \sin\vartheta \vec{L}\vec{L} \, , &
\end{flalign}
\end{subequations}
\endgroup
where
\begingroup
\allowdisplaybreaks
\begin{flalign}
\frac{\partial \vartheta}{\partial \vartheta_I} &= B^I \, , \nonumber&\\
\begin{split}
\frac{\partial \vec{\Lambda}_{1,1} }{\partial \vartheta_I} &=  \cos\vartheta \vec{L}_{,1}B^I  - \sin\vartheta \vec{L} \vartheta_{,1}  B^I + \cos\vartheta \vec{L} B_{,1}^{I}  \\ & + \sin\vartheta (\vec{L} \vec{L}_{,1} +  \vec{L}_{,1} \vec{L} )  B^I + \cos\vartheta \vec{L} \vec{L} \vartheta_{,1} B^I + \sin\vartheta \vec{L} \vec{L} \vartheta_{,1} B_{,1}^I \, .
\end{split} \nonumber&
\end{flalign}
\endgroup
The second derivatives of the rotation matrix $\vec{\Lambda}_1$ and its derivative $\vec{\Lambda}_{1,1}$ with respect to the nodal twist $\vartheta_I$ are given by
\begingroup
\allowdisplaybreaks
\begin{subequations}
\begin{flalign}
\frac{\partial^2 \vec{\Lambda}_1 }{\partial \vartheta_I \partial \vartheta_J} &=  (-\sin\vartheta \vec{L} + \cos\vartheta \vec{L} \vec{L})B^I B^J \, , &\\
\begin{split}
\frac{\partial^2 \vec{\Lambda}_{1,1} }{\partial \vartheta_I \partial \vartheta_J} &=  -\sin\vartheta \vec{L}_{,1}B^I B^J - \cos\vartheta \vec{L} \vartheta_{,1} B^I B^J  - \sin\vartheta \vec L \vartheta_{,1}  B^I B_{,1}^J + (\vec{L} \vec{L}_{,1} +  \vec{L}_{,1} \vec{L} ) \cos\vartheta B^I B^J  
\\ & - \sin\vartheta \vec{L} \vec{L} \vartheta_{,1} B^I B^J + \cos\vartheta \vec{L} \vec{L} B^I B_{,1}^J + \cos\vartheta \vec{L} \vec{L} B_{,1}^I B^J \, .
\end{split} \\ \nonumber &
\end{flalign}
\end{subequations}
\endgroup
\subsection{Internal force derivatives}\label{App:HessianComp}
We derive next the derivatives of the internal forces needed in computing the stiffness matrix. The derivative of the internal force vector $\vec f_{\vec u}$ with respect to the nodal displacements $\vec u_J$ and the twist $\vartheta_J$ are given by
\begingroup
\allowdisplaybreaks
\begin{subequations}
\begin{flalign}
{\dfrac{\partial (\vec f_{\vec u})_I}{\partial \vec u_J}} &= \frac{\partial n}{\partial \vec u_J}\frac{\partial \alpha}{\partial \vec u_I} + n\frac{\partial^2 \alpha}{\partial \vec{u}_I\partial \vec{u}_J} + \sum_{j=2}^3 \left( \frac{\partial m_{j}}{\partial \vec u_J}\frac{\partial \beta_{j}}{\partial \vec u_I} + m_{{j}}\frac{\partial^2 \beta_{{j}}}{\partial \vec{u}_I \partial \vec{u}_J} \right) + \frac{\partial q}{\partial \vec u_J}\frac{\partial \gamma}{\partial \vec u_I} + q\frac{\partial^2 \gamma}{\partial \vec{u}_I \partial \vec{u}_J} \, , \label{eq:Hessian1x1} & \\
{\dfrac{\partial (\vec f_{\vec u})_I}{\partial \vartheta_J}} &=
\sum_{j=2}^3\left( \frac{\partial m_{j}}{\partial \vartheta_J}\frac{\partial \beta_{j}}{\partial \vec u_I} + m_{{j}}\frac{\partial^2 \beta_{{j}}}{\partial \vec{u}_I \partial \vartheta_J}  \right) + \frac{\partial q}{\partial \vartheta_J}\frac{\partial \gamma}{\partial \vec u_I} + q\frac{\partial^2 \gamma}{\partial \vec{u}_I \partial \vartheta_J} \, . \label{eq:Hessian1x2} &
\end{flalign}
\end{subequations}
\endgroup
Due to the symmetry of the second derivatives we have 
\begin{flalign}\label{eq:Hessian2x1}
\dfrac{\partial (f_{\vartheta})_I}{\partial \vec u_J} &= \left[ \dfrac{\partial (\vec f_{\vec u})_I}{\partial \vartheta_J} \right]^\trans \, , &
\end{flalign}
and the derivative of the internal force vector $\vec f_{\vartheta}$ with respect to the nodal twist $\vartheta_J$ is simplified as
\begin{flalign}\label{eq:Hessian2x2}
{\dfrac{\partial (f_{\vartheta})_I}{\partial \vartheta_J}}
&= \sum_{j=2}^3 \left( \frac{\partial m_{j}}{\partial \vartheta_J}\frac{\partial \beta_{j}}{\partial \vartheta_I} + m_{{j}}\frac{\partial^2 \beta_{{j}}}{\partial \vartheta_I \partial \vartheta_J} \right) + \frac{\partial q}{\partial \vartheta_J}\frac{\partial \gamma}{\partial \vartheta_I} + q\frac{\partial^2 \gamma}{\partial \vartheta_I \partial \vartheta_J} \ .&
\end{flalign}
The second derivatives of the strains with respect to the nodal displacements $\vec u_I$ are given by
\begingroup
\allowdisplaybreaks
\begin{subequations}
\begin{flalign}
\frac{\partial^2 \alpha}{\partial \vec{u}_I\partial \vec{u}_J} &= \frac{\partial}{\partial \vec u_J}\left( \vec a_1 \cdot \frac{\partial \vec a_1}{\partial \vec u_I} \right) = \frac{\partial \vec a_1}{\partial \vec u_J} \cdot \frac{\partial \vec a_1}{\partial \vec u_I} \label{eqApp:2ndDeralpha} \, , &\\
\frac{\partial^2 \beta_{\iota}}{\partial \vec{u}_I\partial \vec{u}_J} &= -\frac{\partial}{\partial \vec u_J}\left(\frac{\partial \beta_{\iota}}{\partial \vec u_I} \right) =-\frac{\partial}{\partial \vec u_J}\left( \frac{\partial \vec{a}_{\iota}}{\partial  \vec u_I} \cdot \vec{a}_{1,1} + \vec{a}_{\iota} \cdot \frac{\partial \vec{a}_{1,1}}{\partial \vec u_I} \right) \nonumber \\ &= 
 -\frac{\partial \vec a_{\iota}}{\partial \vec u_I} \cdot \frac{\partial \vec a_{1,1}}{\partial \vec u_J} - \frac{\partial^2 \vec a_{\iota}}{\partial \vec u_I \partial \vec u_J} \cdot \vec a_{1,1} -\frac{\partial \vec a_{\iota}}{\partial \vec u_J} \cdot \frac{\partial \vec a_{1,1}}{\partial \vec u_I}   \label{eqApp:2ndDerbeta} \, , &\\
\frac{\partial^2 \vec a_{\iota}}{\partial \vec u_I \partial \vec u_J} &=\frac{\partial}{\partial \vec u_J}\left( \frac{\partial \vec{\Lambda}_2}{\partial \vec u_I}\vec{\Lambda}_1{\vec A}_{\iota}  \right) = \frac{\partial^2 \vec{\Lambda}_2}{\partial \vec u_I \partial \vec u_J} \vec{\Lambda}_1 \vec A_{\iota} \, , &\\
\frac{\partial^2 \gamma}{\partial \vec{u}_I\partial \vec{u}_J} &= \frac{1}{2} \frac{\partial}{\partial \vec u_J}  \left( \frac{\partial \vec{a}_2}{\partial  \vec u_I} \cdot \vec{a}_{3,1} +
\vec{a}_2 \cdot \frac{\partial \vec{a}_{3,1}}{\partial \vec u_I} \right) \label{eqApp:2ndDergamma} \nonumber \, , &\\
 &= 
\frac{1}{2} \left( \frac{\partial \vec a_2}{\partial \vec u_I} \cdot \frac{\partial \vec a_{3,1}}{\partial \vec u_J} + \frac{\partial^2 \vec a_2}{\partial \vec u_I \partial \vec u_J} \cdot \vec a_{3,1} + \frac{\partial \vec a_2}{\partial \vec u_J} \cdot \frac{\partial \vec a_{3,1}}{\partial \vec u_I} +  \vec a_{2} \cdot \frac{\partial^2 \vec a_{3,1}}{\partial \vec u_I \partial \vec u_J} \right) \, , &\\
\frac{\partial^2 \vec a_{3,1}}{\partial \vec u_I \partial \vec u_J} &=\frac{\partial}{\partial \vec u_J}\left( 
\frac{\partial \vec{\Lambda}_{2,1}}{\partial \vec u_I}\vec{\Lambda}_1{\vec A}_3 +
\frac{\partial \vec{\Lambda}_2}{\partial \vec u_I}\vec{\Lambda}_{1,1}{\vec A}_3 + \frac{\partial \vec{\Lambda}_2}{\partial \vec u_I}\vec{\Lambda}_1{\vec A}_{3,1}
   \right) \nonumber \, , &\\
&= \frac{\partial^2 \vec{\Lambda}_{2,1}}{\partial \vec u_I \partial \vec u_J} \vec{\Lambda}_1 \vec A_3 + \frac{\partial^2 \vec{\Lambda}_2}{\partial \vec u_I \partial \vec u_J} \vec{\Lambda}_{1,1} \vec A_3 + \frac{\partial^2 \vec{\Lambda}_2}{\partial \vec u_I \partial \vec u_J} \vec{\Lambda}_1 \vec A_{3,1} \, .  &
\end{flalign}
\end{subequations}
\endgroup
And the mixed second derivatives of the strains are given by
\begingroup
\allowdisplaybreaks
\begin{subequations}
\begin{flalign}
\frac{\partial^2 \beta_{\iota}}{\partial \vec{u}_I\partial \vartheta_J} &= -\frac{\partial}{\partial \vartheta_J}\left(\frac{\partial \beta_{\iota}}{\partial \vec u_I} \right) =-\frac{\partial}{\partial \vartheta_J}\left( \frac{\partial \vec{a}_{\iota}}{\partial  \vec u_I} \cdot \vec{a}_{1,1} + \vec{a}_{\iota} \cdot \frac{\partial \vec{a}_{1,1}}{\partial \vec u_I} \right) \nonumber &\\ &= 
-\frac{\partial^2 \vec a_{\iota}}{\partial \vec u_I \partial \vartheta_J} \cdot \vec a_{1,1} -\frac{\partial \vec a_{\iota}}{\partial \vartheta_J} \cdot \frac{\partial \vec a_{1,1}}{\partial \vec u_I} \, , &\\
\frac{\partial^2 \vec a_{\iota}}{\partial \vec u_I \partial \vartheta_J} &=\frac{\partial}{\partial \vartheta_J}\left( \frac{\partial \vec{\Lambda}_2}{\partial \vec u_I}\vec{\Lambda}_1{\vec A}_{\iota}  \right) = 
\frac{\partial \vec{\Lambda}_2}{\partial \vec u_I} \frac{\partial \vec{\Lambda}_1}{\partial \vartheta_J} \vec A_{\iota} \, , &\\
\frac{\partial^2 \gamma}{\partial \vec{u}_I\partial \vartheta_J} &= \frac{1}{2} \frac{\partial}{\partial \vartheta_J}  \left( \frac{\partial \vec{a}_2}{\partial  \vec u_I} \cdot \vec{a}_{3,1} +
\vec{a}_2 \cdot \frac{\partial \vec{a}_{3,1}}{\partial \vec u_I} \right) \nonumber &\\
 &= 
\frac{1}{2} \left( \frac{\partial \vec a_2}{\partial \vec u_I} \cdot \frac{\partial \vec a_{3,1}}{\partial \vartheta_J} + \frac{\partial^2 \vec a_2}{\partial \vec u_I \partial \vartheta_J} \cdot \vec a_{3,1} + \frac{\partial \vec a_2}{\partial\vartheta_J} \cdot \frac{\partial \vec a_{3,1}}{\partial \vec u_I} +  \vec a_{2} \cdot \frac{\partial^2 \vec a_{3,1}}{\partial \vec u_I \partial \vartheta_J} \right) \, , &\\
\frac{\partial^2 \vec a_{3,1}}{\partial \vec u_I \partial \vartheta_J} &=\frac{\partial}{\partial \vartheta_J}\left( 
\frac{\partial\vec{\Lambda}_{2,1}}{\partial \vec u_I}\vec{\Lambda}_1{\vec A}_3 +
\frac{\partial \vec{\Lambda}_2}{\partial \vec u_I}\vec{\Lambda}_{1,1}{\vec A}_3 + \frac{\partial \vec{\Lambda}_2}{\partial \vec u_I}\vec{\Lambda}_1{\vec A}_{3,1}
   \right) \nonumber &\\
&= \frac{\partial \vec{\Lambda}_{2,1}}{\partial \vec u_I } \frac{\partial \vec{\Lambda}_1}{\partial \vartheta_J} \vec A_3 +
\frac{\partial \vec{\Lambda}_2}{\partial \vec u_I} \frac{\partial \vec{\Lambda}_{1,1}}{\partial \vartheta_J} \vec A_3 +
\frac{\partial \vec{\Lambda}_2}{\partial \vec u_I} \frac{\partial \vec{\Lambda}_1}{\partial \vartheta_J} \vec A_{3,1} \, . &
\end{flalign}
\end{subequations}
\endgroup
\\
Lastly, the second derivatives of the strains with respect to the nodal twist $\vartheta_I$ take the following forms.
\begingroup
\allowdisplaybreaks
\begin{subequations}
\begin{flalign}
\frac{\partial^2 \beta_{\iota}}{\partial \vartheta_I\partial \vartheta_J} &= -\frac{\partial}{\partial \vartheta_J}\left(\frac{\partial \beta_{\iota}}{\partial \vartheta_I} \right) =-\frac{\partial}{\partial \vartheta_J}\left( \frac{\partial \vec{a}_{\iota}}{\partial  \vartheta_I} \cdot \vec{a}_{1,1}\right) = 
-\frac{\partial^2 \vec a_{\iota}}{\partial \vartheta_I \partial \vartheta_J} \cdot \vec a_{1,1} \, , &\\
\frac{\partial^2 \vec a_{\iota}}{\partial \vartheta_I \partial \vartheta_J} &= \frac{\partial}{\partial \vartheta_J} \left( \vec{\Lambda}_2 \frac{\partial \vec{\Lambda}_1}{\partial \vartheta_I}{\vec A}_{\iota} \right) =  
 \vec{\Lambda}_2 \frac{\partial^2 \vec{\Lambda}_1}{\partial \vartheta_I \partial \vartheta_J}{\vec A}_{\iota} \, , & \\
\frac{\partial^2 \gamma}{\partial \vartheta_I \partial \vartheta_J} &= \frac{1}{2} \frac{\partial}{\partial \vartheta_J}  \left( \frac{\partial \vec{a}_2}{\partial \vartheta_I} \cdot \vec{a}_{3,1} +
\vec{a}_2 \cdot \frac{\partial \vec{a}_{3,1}}{\partial \vartheta_I} \right) \nonumber \\
 &= 
\frac{1}{2} \left( \frac{\partial \vec a_2}{\partial \vartheta_I} \cdot \frac{\partial \vec a_{3,1}}{\partial \vartheta_J} + \frac{\partial^2 \vec a_2}{\partial \vartheta_I \partial \vartheta_J} \cdot \vec a_{3,1} + \frac{\partial \vec a_2}{\partial\vartheta_J} \cdot \frac{\partial \vec a_{3,1}}{\partial \vartheta_I} +  \vec a_{2} \cdot \frac{\partial^2 \vec a_{3,1}}{\partial \vartheta_I \partial \vartheta_J} \right) \, , \\
\frac{\partial^2 \vec a_{3,1}}{\partial \vartheta_I \partial \vartheta_J} &= \frac{\partial}{\partial \vartheta_J}\left( 
\vec{\Lambda}_{2,1}\frac{\partial \vec{\Lambda}_1}{\partial \vartheta_I}{\vec A}_3 + 
\vec{\Lambda}_2 \frac{\partial \vec{\Lambda}_{1,1}}{\partial \vartheta_I}{\vec A}_3 + 
\vec{\Lambda}_2\frac{\partial \vec{\Lambda}_1}{\partial \vartheta_I}{\vec A}_{3,1}
   \right) \nonumber \\
&= \vec{\Lambda}_{2,1}\frac{\partial^2 \vec{\Lambda}_1}{\partial \vartheta_I \partial \vartheta_J}{\vec A}_3 + 
\vec{\Lambda}_2 \frac{\partial^2 \vec{\Lambda}_{1,1}}{\partial \vartheta_I \partial \vartheta_J}{\vec A}_3 + 
\vec{\Lambda}_2\frac{\partial^2 \vec{\Lambda}_1}{\partial \vartheta_I \partial \vartheta_J}{\vec A}_{3,1} \, .
\end{flalign}
\end{subequations}
\endgroup
\section{Geometry definition of a weft-knitted RVE}\label{app:KnitGeometry}
\begin{figure} 
\includegraphics[scale=0.55]{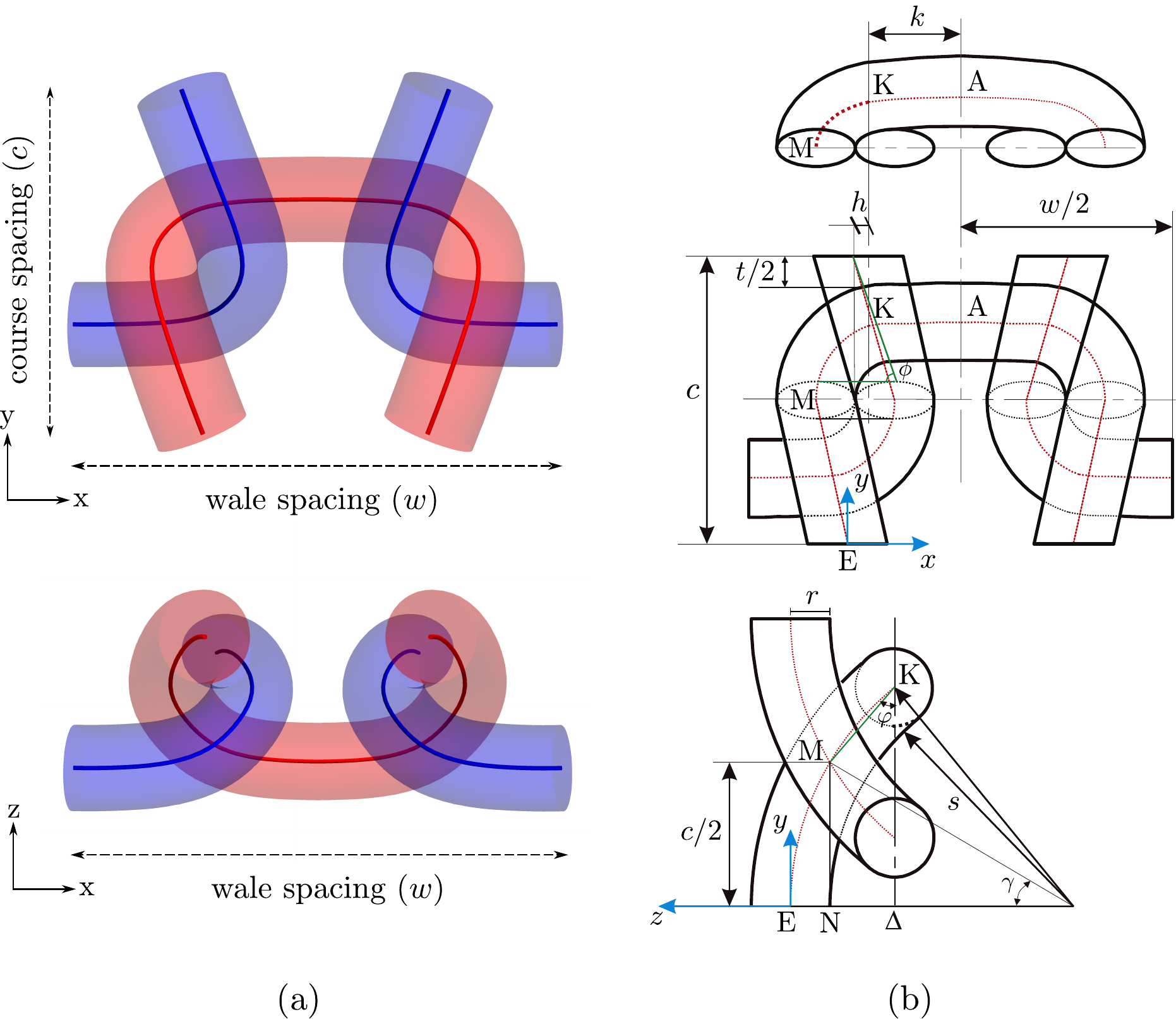}
\centering
\caption[Different views of a RVE of a weft-knitted textile]{Geometry modelling of a weft-knitted RVE (a) Different views of a RVE of a weft-knitted textile (b) Introducing the geometric parameters that define the geometric description of the selected RVE~\cite{vassiliadis2007}.}
\label{fig:knitRVEgeometry}
\end{figure}

The following equations for the geometry depicted in Figure~\ref{fig:knitRVEgeometry} are extracted from Vassiliadis et al.~\cite{vassiliadis2007a} and smoothened at $C^0$ locations to recover the real smooth geometry of the RVE.

In part EM $(0 < y < c/2)$:
\begin{subequations}
\begin{align}
x(y) &= -\dfrac{2r}{c}y \, ,\\
z(y) &= \sqrt{\left(s + r \right)^2 - y^2} - \left(s + r\right) \, ,
\end{align}
\end{subequations}
where $s = \dfrac{1}{4r}\left[ \left( c - r - \dfrac{t}{2}\right)^2 - \left( r + \dfrac{t}{2}\right)^2 \right]$.

In part MK $( c/2 < y < c/2 + R)$:
\begin{subequations}
\begin{align}
x(y) &= h - a \sqrt{1 - \left( \dfrac{y-c/2}{b} \right)} \, , \\
z(y) &= \sqrt{\left(s + r \right)^2 - y^2} - \left(s + r\right) \, ,
\end{align}
\end{subequations}
where \mbox{$h = \left( \dfrac{c}{2} - R \right)\tan{\left( \dfrac{\pi}{2} - \phi \right)}$}, \mbox{$R = \dfrac{c}{2} - \dfrac{t}{2} - r$}, \mbox{$\phi = \arctan \left( \dfrac{c-2r\sin{\gamma}}{2r} \right)$} and \mbox{$\gamma = \arcsin{ \left( \dfrac{c/2}{s + r} \right) }$}.

In part KA $( x(y = c/2 + R) < x < w/4)$:
\begin{subequations}
\begin{align}
y(z) &= \sqrt{\left(s + r \right)^2 - \left(z + s + r \right)^2} \, , \\
z(x) &= OZ - \sqrt{A^2 - (x - OX)^2} \, ,
\end{align}
\end{subequations}
where 
\[ 
OZ = \left[(x_2 - OX)^2 - (x_1 - OX)^2 + z^2_2 - z^2_1\right]/2(z_2 - z_1) \, , 
\] 
 \[ A = \sqrt{(x_1 - OX)^2 + (z_1 - OZ)^2} \, , \] 
 \[ OX = w/4 \, .\]
 Additionally, $(x_1, z_1)$ and $(x_2, z_2)$ are the coordinates of two points in part MK with \mbox{$y_1 = c/2 + R - 0.001$} and \mbox{$y_2 = c/2 + R$}.

A cubic B-spline curve was fitted at ${C}^0$ continuous locations namely, M and K. Moreover, this smooth curve EMKA is mirrored on the vertical axis passing at point A to obtain the complete geometry of the red yarn shown in Figure \ref{fig:knitRVEgeometry}. Lastly, one smooth EMKA is rotated and translated appropriately to create the blue left yarn of Figure \ref{fig:knitRVEgeometry} and mirrored subsequently to complete the RVE geometry.

\section*{Data availability statement}
The data that support the findings of this study are available from the corresponding author upon reasonable request.

\bibliographystyle{unsrt}
\bibliography{knitting}

\end{document}